\documentclass[a4paper,12pt,onehalfspace]{article}
\usepackage{epsfig,graphicx}
\usepackage{psfrag}
\usepackage{amssymb,bbm}
\usepackage{amsmath}
\usepackage{amsfonts}
\usepackage{setspace}
\usepackage[toc,page]{appendix}
\usepackage{enumerate}
\usepackage{tikz}
\usetikzlibrary{patterns,shapes,snakes,backgrounds,calc}
\tikzset{axis line style/.style={thin, gray, -stealth}}
\newcommand*{\TickSize}{2pt}%
\usepackage{float}
\usepackage{tikz-cd}
\usepackage{blindtext}
\usepackage{upgreek}
\usepackage{fancyhdr}
\usepackage{subcaption}

\usepackage{pstricks,pst-math,pst-plot}
\usepackage[dotinlabels]{titletoc}
\usepackage{psfrag}
\usepackage{pstricks-add}

\fancyheadoffset{0cm}

\fancypagestyle{plain}{%
\fancyhf{}%
\fancyhead[R]{\thepage}%
}

\newcommand*{\xMin}{-6}%
\newcommand*{\xMax}{7}%
\newcommand*{\yMin}{-5}%
\newcommand*{\yMax}{4}%
\usepackage[colorlinks]{hyperref}
\newcounter{theorem}
\newtheorem{theorem}{Theorem}
\newtheorem{example}{Example}
\newtheorem{remark}{Remark}
\newtheorem{proposition}{Proposition}
\newtheorem{lemma}{Lemma}
\newtheorem{definition}{Definition}

\newenvironment{proof}[1][Proof]{\textbf{#1.} }{\rule{0.5em}{0.5em}}

\title{Trimming of Finite Subsets of the Manhattan Plane\date{}}
\author{G\"{o}k\c{c}e \c{C}AKMAK \footnote{Eskisehir Technical University, Science Faculty, Department of Mathematics, 26470, Eski\c{s}ehir, Turkey.  e-mail: gokcecakmak@eskisehir.edu.tr}
\thanks{Corresponding Author.}
\and Ali DEN\.{I}Z \footnote{Eskisehir Technical University, Science Faculty, Department of Mathematics, 26470, Eski\c{s}ehir, Turkey. e-mail: adeniz@eskisehir.edu.tr } \and  \c{S}ahin KO\c{C}AK \footnote{Anadolu University (emeritus), Science Faculty, Department of Mathematics, 26470, Eski\c{s}ehir, Turkey. e-mail: skocak@anadolu.edu.tr}}
\begin{document}
\maketitle
\thispagestyle{empty}
%%%%%%%%%%

%\textbf{MSC2010:}
\begin{abstract}
V. Turaev defined recently an operation of ``Trimming" for pseudo-metric spaces and analysed the tight span of (pseudo-)metric spaces via this process. In this work we investigate the trimming of finite subspaces of the Manhattan plane. We show that this operation amounts for them to taking the metric center set and we give an algorithm to construct the tight spans via trimming.
\end{abstract}
%\maketitle

\section{Introduction}
Our aim in this work is to study the finite subsets of the Manhattan plane from the viewpoint of trimming theory introduced recently by V. Turaev \cite{Turaev2016,Turaev2018}. We want first to explain briefly the notions of trimming, trimming sequence and trimming cylinder defined by Turaev.

Let $(X,d)$ be a non-empty finite metric space. (We use the notion of metric space in the genuine sense that $d(x,y)=0$ if and only if $x=y$. We sometimes simply write $X$ for $(X,d)$). The process of trimming introduced by V. Turaev (\cite{Turaev2018}) can be explained as follows. First assume $|X|\geq 3$. We define an equivalence relation $\Re$ on $X$:
 \[
 x\Re y \; \mbox{for} \; x, y \in X :\Leftrightarrow\;  x=y \quad \text{or} \quad d(x,y)=\underline{d}(x)+\underline{d}(y)
 \]
whereby
 \[
       \underline{d}(x)=\min_{\substack{y,z\in X\setminus\{x\}\\ y\neq z}} \dfrac{d(x,y)+d(x,z)-d(y,z)}{2}.
       \]

 We then define the following metric $d^1$ on the set $X/\Re $ of equivalence classes by $d^1(\overline{x},\overline{y})=d(x,y)-\underline{d}(x)-\underline{d}(y)$ for $\overline{x},\overline{y}\in X/\Re$. The metric space $(X/\Re, d^1)$ is said to be obtained by trimming of $(X,d)$ and is denoted by $(t(X),d^1)$ (or simply by $t(X)$). The map $p:(X,d)\to (t(X),d^1), p(x)=\overline{x}$, is called the trimming projection. For $|X|=1$ or $|X|=2$ we define $t(X)$ to be a one-point space, the trimming projection being the trivial one.
\paragraph{Some Remarks:}
\begin{enumerate}
  \item For convenience, one can define $\underline{d}(x)=\underline{d}(y)=\dfrac{1}{2}d(x,y)$ for $X=\{x, y\}$ and $\underline{d}(x)=0$ for $X=\{x\}$.
  \item The quantity $\triangle_{x,yz}=\dfrac{1}{2}(d(x,y)+d(x,z)-d(y,z))$ for $x,y,z \in X$ is called the Gromov product at $x$ with respect to $y$ and $z$. For a discussion on Gromov products see \cite{Bilge2017}.
  \item The process of trimming can be defined for any pseudo-metric space (\cite{Turaev2018}).
  \item $\underline{d}(x)$ is called the pendant length of $x\in X$. If $\underline{d}(x)=0$ for all $x\in X$, then $(X,d)$ is called pendant-free or trim.

      For distinct $x,y,z \in X$, $x$ is called ``between" $y$ and $z$ in the Menger sense, if $d(y,z)=d(y,x)+d(x,z)$. The pendant length $\underline{d}(x)$ vanishes if and only if $x$ is between two other points $y$ and $z$ (for $|X|\geq 3$).
      \item For a 3-point space $X=\{x,y,z\}$, $t(X)$ is a singleton: Since $\underline{d}(x)=\triangle_{x,yz}, \underline{d}(y)=\triangle_{y,xz}$ and $\underline{d}(z)=\triangle_{z,xy}$, $x\Re y$ because $d(x,y)=\underline{d}(x)+\underline{d}(y)$; respectively, $x\Re z$ and $y\Re z$.
\end{enumerate}

One should not get the impression that the process of trimming produces a trim (pendant-free) space, as the following example shows.
\begin{example}
Let $X=\{1,2,3,4,5\}\subset \mathbb{R}$ with the metric induced from the standard metric of $\mathbb{R}$. Then $t(X)=\{\overline{1}, \overline{3}, \overline{4}\}$ and $\underline{d}^1(\overline{1})=1$, so that $(t(X),d^1)$ is not trim (see Figure \ref{fig1}).
\begin{figure}[h]
\centering
\begin{tikzpicture}
\draw [fill] (0,0) circle (.05); \node [below left] at (0,0) {$\overline{1}=\{1,2\}$};
\draw [fill] (1,0) circle (.05); \node [below] at (1,0) {$\overline{3}=\{3\}$}; \node at (-2.5,0) {$t(X)$};
\draw [fill] (2,0) circle (.05); \node [below right] at (2,0) {$\overline{4}=\{4,5\}$};
\draw [fill] (-1,1.5) circle (.05); \node [above left] at (-1,1.5) {1};
\draw [fill] (0,1.5) circle (.05); \node [above] at (0,1.5) {2}; \node at (-2.5,1.5) {$X$};
\draw [fill] (1,1.5) circle (.05);\node [above] at (1,1.5) {3};
\draw [fill] (2,1.5) circle (.05);\node [above] at (2,1.5) {4};
\draw [fill] (3,1.5) circle (.05);\node [above] at (3,1.5) {5};
 \draw [->] (-1,1.3) -- (-0.3,0.2);
 \draw [->]  (0,1.3) -- (0,0.2);
  \draw[->] (2,1.3) -- (2,0.2);
 \draw [->]  (3,1.3) -- (2.2,0.2);
  \draw [->] (1,1.5) -- (1,0.2);
\end{tikzpicture}
 \caption{A representation of $(X,d)$ and $(t(X),d^1)$.}\label{fig1}
 \end{figure}
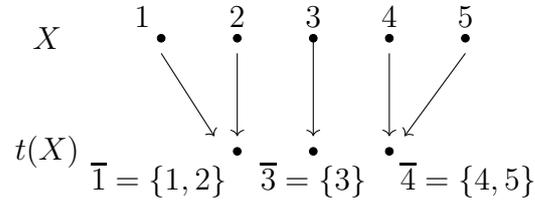
\end{example}

The goal of the process of trimming is to obtain a trim (pendant-free) space by applying this operation succesively. The consecutive trimmings of $(X,d)$ create metric spaces $(t^i(X),d^i) \;(i>0)$ and non-expansive surjections $p_i: t^i(X)\to t^{i+1}(X), p_i(x_{i})=x_{i+1}$, so that we get a sequence
\[
X=t^0(X)\xrightarrow{p_0=p}t(X)\xrightarrow{p_1}t^2(X)\rightarrow \cdots \rightarrow t^i(X)\xrightarrow{p_i}t^{i+1}(X)\rightarrow \cdots
\]
which is called the trimming sequence of $X$. Starting from a point $x=x_{0}\in X$, we get a sequence $(x_{i})_{i\geq0}$, $x_i\in t^i(X)$ with $x_{i+1}=p_i(x_{i})$, which is called the trimming sequence of $x$. We write sometimes $X_i$ for $t^i(X)$. For $x_{i},y_{i}\in X_i, \; (i>0)$ we have   $d^i(x_{i},y_{i})=d^{i-1}(x_{i-1},y_{i-1})-\underline{d}^{i-1}(x_{i-1})-\underline{d}^{i-1}(y_{i-1})$.

If $X$ is trim, then $t(X)=X$ (or naturally isometric if one so wishes); thus, starting from an arbitrary $X$, if one arrives at a trim $t^i(X)$, then the trimming sequence stabilizes from there on.

\begin{figure}[h]
\begin{center}
\scalebox{1}{\begin{tikzpicture}
\draw [fill] (0,0) circle (.05); \node [below left] at (0,0) {$x$};
\draw [fill] (4,0) circle (.05); \node [below right] at (4,0) {$y$};
\draw [fill] (4,3) circle (.05); \node [above right] at (4,3) {$z$};
\draw [fill] (0,3) circle (.05); \node [above left] at (0,3) {$t$};

\node [above left] at (.5,.5) {1};
\node [below] at (2,1) {2};
\node [above right] at (3.5,.5) {1};
\node [right] at (3,1.5) {1};
\node [below right] at (3.5,2.5) {1};
\node [above] at (2,2) {2};
\node [below left] at (.5,2.5) {1};
\node [left] at (1,1.5) {1};

\node [below] at (2,-1) {$X=\{x,y,z,t\}$};

\draw (1,1) circle (.05);
\draw (3,1) circle (.05);
\draw (3,2) circle (.05);
\draw (1,2) circle (.05);

 \draw (0,0) -- (1,1);
 \draw (1,1) -- (3,1);
  \draw (3,1) -- (4,0);
 \draw (3,1) -- (3,2);
  \draw (3,2) -- (4,3);
  \draw (3,2) -- (1,2);
  \draw (1,2) -- (0,3);
  \draw (1,2) -- (1,1);

\end{tikzpicture}}
\qquad \qquad \scalebox{1}{
\begin{tikzpicture}
\draw [fill] (1,1) circle (.05); \node [below left] at (1,1) {$p(x)$};
\draw [fill] (3,1) circle (.05); \node [below right] at (3,1) {$p(y)$};
\draw [fill] (3,2) circle (.05); \node [above right] at (3,2) {$p(z)$};
\draw [fill] (1,2) circle (.05); \node [above left] at (1,2) {$p(t)$};

\node [below] at (2,1) {2};
\node [right] at (3,1.5) {1};
\node [above] at (2,2) {2};
\node [left] at (1,1.5) {1};

\node [below] at (2,-1) {$t(X)=\{p(x),p(y),p(z),p(t)\}$};

 \draw (1,1) -- (3,1);
 \draw (3,1) -- (3,2);
  \draw (3,2) -- (1,2);
  \draw (1,2) -- (1,1);
\end{tikzpicture}}
\end{center}
\caption{A bijective trimming projection $p:X\to t(X)$ with non-trim $X$ and trim $t(X)$.}\label{fig2}
 \end{figure}
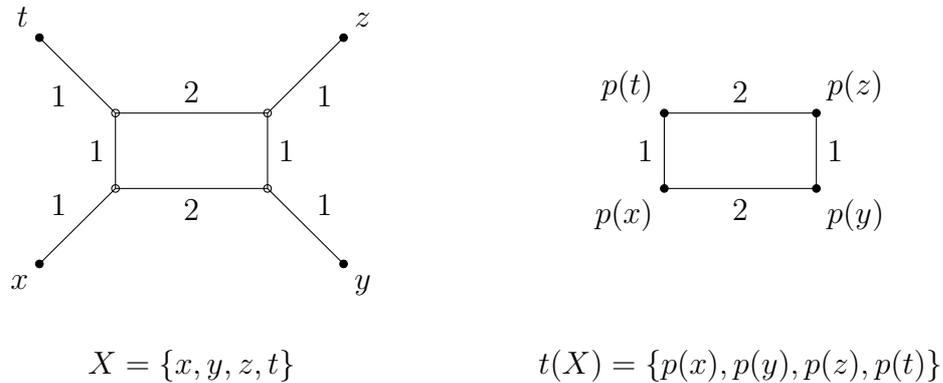

If the trimming projection $p:X\to t(X)$ is bijective then $t(X)$ is trim (though $X$ itself might not be). As an example see Figure \ref{fig2}, where the distances are to be read as the lengths of shortest paths on the given auxiliary graphs.

%%% 그래픽 추가 캡션을 다세요
By this reason, successive trimming of a finite metric space will eventually produce a trim space since the sequence $|t^i(X)|_{i\geq 0}$ is non-increasing and thus a bijective $p_i$ will emerge.

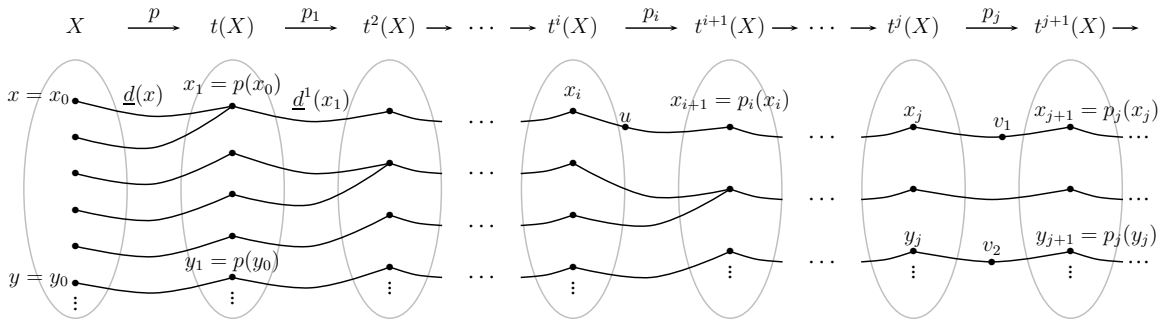
\begin{figure}[h]
\resizebox{\linewidth}{!}{
\begin{pspicture}(-1,4)(20.5,-3.5)
%\psset{unit=0.65cm}
%\psgrid
\psellipse[linecolor=lightgray](0,0)(1,2.5)
\psellipse[linecolor=lightgray](3,0)(1,2.5)
\psellipse[linecolor=lightgray](6,0)(1,2.5)
\psellipse[linecolor=lightgray](9.5,0)(1,2.5)
\psellipse[linecolor=lightgray](12.5,0)(1,2.5)
\psellipse[linecolor=lightgray](16,0)(1,2.5)
\psellipse[linecolor=lightgray](19,0)(1,2.5)

 %\psellipse[linecolor=gray](0,0)(1,2.5)
 \psdots*(0,1.7)(0,1)(0,0.3)(0,-0.4)(0,-1.1)(0,-1.8)%kümedeki noktalar
 \psdots*[dotscale=0.4](0,-2.1)(0,-2.2)(0,-2.3)%kümenin altındaki 3 küçük nokta
 \pscurve(0,1.7)(1.5,1.4)(3,1.6) %noktaları birleştiren eğriler
 \pscurve(0,1)(1.5,0.8)(3,1.6)
 \pscurve(0,0.3)(1.5,0.1)(3,0.7)
 \pscurve(0,-0.4)(1.5,-0.6)(3,-0.1)
 \pscurve(0,-1.1)(1.5,-1.3)(3,-0.9)
 \pscurve(0,-1.8)(1.5,-2)(3,-1.7)
 \rput[B](0,3){$X$}
 \psline{->}(1,3.1)(2,3.1)
\rput[B](1.5,3.3){$p$}
\rput[B](-0.7,1.7){$x=x_0$} \rput[B](-0.7,-1.8){$y=y_0$}

 %\psellipse[linecolor=gray](3,0)(1,2.5)
 \psdots*(3,1.6)(3,0.7)(3,-0.1)(3,-0.9)(3,-1.7)
 \psdots*[dotscale=0.4](3,-2)(3,-2.1)(3,-2.2)
 \pscurve(3,1.6)(4.5,1.3)(6,1.5)
 \pscurve(3,0.7)(4.5,0.3)(6,0.5)
 \pscurve(3,-0.1)(4.5,-0.3)(6,0.5)
 \pscurve(3,-0.9)(4.5,-1.1)(6,-0.5)
 \pscurve(3,-1.7)(4.5,-1.9)(6,-1.5)
\rput[B](3,3){$t(X)$}
 \psline{->}(4,3.1)(5,3.1)
\rput[B](4.5,3.3){$p_1$}
\rput[B](3,1.9){$x_1=p(x_0)$}  \rput[B](1.3,1.7){$\underline{d}(x)$} \rput[B](4.7,1.6){$\underline{d}^1(x_1)$}
\rput[B](3,-1.5){$y_1=p(y_0)$}

 %\psellipse[linecolor=gray](6,0)(1,2.5)
 \psdots*(6,1.5)(6,0.5)(6,-0.5)(6,-1.5)
 \psdots*[dotscale=0.4](6,-1.8)(6,-1.9)(6,-2)
 \rput[B](6,3){$t^2(X)$}
 %first arms
 \pscurve(6,1.5)(6.5,1.35)(7,1.3)
 \pscurve(6,0.5)(6.5,0.35)(7,0.3)
 \pscurve(6,-0.5)(6.5,-0.65)(7,-0.7)
 \pscurve(6,-1.5)(6.5,-1.65)(7,-1.7)
  \psline{->}(6.7,3.1)(7.2,3.1)
 %2nd arms
 \pscurve(8.5,1.3)(9,1.35)(9.5,1.5)
 \pscurve(8.5,0.3)(9,0.35)(9.5,0.5)
 \pscurve(8.5,-0.7)(9,-0.65)(9.5,-0.5)
 \pscurve(8.5,-1.7)(9,-1.65)(9.5,-1.5)
  \psline{->}(8.3,3.1)(8.8,3.1)
 %dots between arms
 \psdots*[dotscale=0.4](7.55,3.1)(7.75,3.1)(7.95,3.1)(7.55,1.3)(7.75,1.3)(7.95,1.3)(7.55,0.3)(7.75,0.3)(7.95,0.3)(7.55,-0.7)(7.75,-0.7)(7.95,-0.7)(7.55,-1.7)(7.75,-1.7)(7.95,-1.7)

 %\psellipse[linecolor=gray](9.5,0)(1,2.5)
\psdots*(9.5,1.5)(9.5,0.5)(9.5,-0.5)(9.5,-1.5)(10.5,1.2)
 \psdots*[dotscale=0.4](9.5,-1.8)(9.5,-1.9)(9.5,-2)
 \rput[B](10.5,1.3){$u$}
 \pscurve(9.5,1.5)(11,1.1)(12.5,1.2)
 \pscurve(9.5,0.5)(11,-0.15)(12.5,0)
 \pscurve(9.5,-0.5)(11,-0.7)(12.5,0)
 \pscurve(9.5,-1.5)(11,-1.7)(12.5,-1.2)
  \rput[B](9.5,3){$t^i(X)$}
  \psline{->}(10.5,3.1)(11.5,3.1)
  \rput[B](11,3.3){$p_i$}
\rput[B](9.5,1.8){$x_i$}

 %\psellipse[linecolor=gray](12.5,0)(1,2.5)
\psdots*(12.5,1.2)(12.5,0)(12.5,-1.2)
 \psdots*[dotscale=0.4](12.5,-1.5)(12.5,-1.6)(12.5,-1.7)
  \rput[B](12.5,3){$t^{i+1}(X)$}
  \rput[B](12.5,1.6){$x_{i+1}=p_i(x_i)$}
  %first arms
 \pscurve(12.5,1.2)(13,1.05)(13.5,1)
 \pscurve(12.5,0)(13,-0.15)(13.5,-0.2)
 \pscurve(12.5,-1.2)(13,-1.35)(13.5,-1.4)
  \psline{->}(13.3,3.1)(13.8,3.1)
 %2nd arms
 \pscurve(15,1)(15.5,1.05)(16,1.2)
 \pscurve(15,-0.2)(15.5,-0.15)(16,0)
 \pscurve(15,-1.4)(15.5,-1.35)(16,-1.2)
  \psline{->}(14.8,3.1)(15.3,3.1)
 %dots between arms
 \psdots*[dotscale=0.4](14.05,3.1)(14.25,3.1)(14.45,3.1)(14.05,1)(14.25,1)(14.45,1)(14.05,-0.2)(14.25,-0.2)(14.45,-0.2)(14.05,-1.4)(14.25,-1.4)(14.45,-1.4)

 %\psellipse[linecolor=gray](16,0)(1,2.5)
\psdots*(16,1.2)(16,0)(16,-1.2)(17.5,-1.4)(17.7,1)
\psdots*[dotscale=0.4](16,-1.5)(16,-1.6)(16,-1.7)
\pscurve(16,1.2)(17.5,1)(19,1.2)
\pscurve(16,0)(17.5,-0.2)(19,0)
\pscurve(16,-1.2)(17.5,-1.4)(19,-1.2)
\rput[B](16,3){$t^j(X)$}
\psline{->}(17,3.1)(18,3.1)
\rput[B](17.5,3.3){$p_j$}
\rput[B](17.7,1.2){$v_1$}
\rput[B](17.5,-1.2){$v_2$}
\rput[B](16,1.4){$x_j$}
\rput[B](16,-1){$y_j$}
 %\psellipse[linecolor=gray](19,0)(1,2.5)
 \psdots*(19,1.2)(19,0)(19,-1.2)
 \psdots*[dotscale=0.4](19,-1.5)(19,-1.6)(19,-1.7)
 \pscurve(19,1.2)(19.5,1.05)(20,1)
 \pscurve(19,0)(19.5,-0.15)(20,-0.2)
 \pscurve(19,-1.2)(19.5,-1.35)(20,-1.4)
 \rput[B](19,3){$t^{j+1}(X)$}
 \rput[B](19.5,1.4){$x_{j+1}=p_j(x_j)$}
 \rput[B](19.5,-1){$y_{j+1}=p_j(y_j)$}
 \psline{->}(19.8,3.1)(20.3,3.1)
%dots after arms (20.31,3.1)(20.35,3.1)(20.39,3.1)
 \psdots*[dotscale=0.4](20.15,1)(20.30,1)(20.45,1)(20.15,-0.2)(20.30,-0.2)(20.45,-0.2)(20.15,-1.4)(20.30,-1.4)(20.45,-1.4)
\end{pspicture}
}
\caption{The trimming cylinder of $(X,d)$ (In the figure, $v_1$ is below $u$ and $v_2$ is not below $u$).}\label{trimmingcylinder}
\end{figure}

Now we want to recall another notion introduced by Turaev (\cite{Turaev2018}). The trimming cylinder $C=C(X)$ of a finite metric space $X$ is a pseudo-metric graph that can be defined as follows: The vertex set is $\bigsqcup_{i\geq0} t^i(X)$ and the edges are put between the vertices $x_{i}\in t^i(X)$ and $p_i(x_{i})\in t^{i+1}(X)$ with $\underline{d}^i(x_{i})$ assigned as weight (see Figure \ref{trimmingcylinder}). The pseudo-metric $\rho$ on $C(X)$ can be defined as follows (for more detail cf. \cite{Turaev2018}):
\begin{align*}
  \rho &: C(X)\times C(X)\to \mathbb{R}\\
  \rho(u,v)&=\begin{cases}
               \displaystyle \sum_{n=i}^j\underline{d}^n(x_n)-\lambda\underline{d}^i(x_i)-(1-\mu)\underline{d}^j(x_j), \quad \mbox{if}\; $v$\; \mbox{is below}\; $u$  \\
               \\
              \displaystyle d_1(x,y)-\sum_{n=0}^{i-1}\underline{d}^n(x_n)-\sum_{n=0}^{j-1}\underline{d}^n(y_n)-\lambda\underline{d}^i(x_i)-\mu\underline{d}^j(y_j), \ \mbox{otherwise}
             \end{cases}
\end{align*}
whereby $u$ is on a path determined by a trajectory $(x_i\in t^i(X))_{i\geq 0}$ (with $x_0=x\in X$) and given by $u=(1-\lambda)x_i+\lambda x_{i+1}$ (thus $u$ lying on the edge connecting $x_i$ with $x_{i+1}$ with weight $\underline{d}^i(x_i)$, $0\leq \lambda \leq 1$) and similarly $v=(1-\mu)y_j+\mu y_{j+1}$. ``$v$ is below $u$" means that $v$ is on a path starting from an $x\in X$ and occuring after $u$ (see Figure \ref{trimmingcylinder}).

We denote the subgraph restricted to $\bigsqcup_{i=0}^{k} t^i(X)$ (with $k>0$) by $C_k(X)$. The metric quotient of the trimming cylinder is denoted as $\overline{C(X)}$.

%% definition of tight span
\par Let us now define the notion of tight span for a finite metric space. Consider the set of functions $f:X\to \mathbb{R}^{\geq 0}$ satisfying the following two conditions:
\begin{itemize}
  \item [i.] $f(x)+f(y)\geq d(x,y)$ for all $x,y\in X$,
  \item [ii.] For each $x\in X$ there exists $y\in X$ such that $f(x)+f(y)= d(x,y)$.
\end{itemize}

The tight span $T(X)$ of $X$ is then this set of functions with the maximum metric
\begin{equation*}
  d_{\infty}(f,g)=\max_{x\in X} |f(x)-g(x)|.
\end{equation*}
If $X$ is a subset of the Manhattan plane, then a closed, geodesically convex subspace $Y\supset X$ of the Manhattan plane, which is minimal with respect to these properties is isometric to the tight span $T(X)$ of $X$ (\cite{Kilic2016}).

In \cite{Turaev2018}, Turaev defined a trimming filtration
\begin{equation*}
    T(X)\supset T(X_1) \supset T(X_2) \supset \cdots \supset T(X_i) \supset T(X_{i+1}) \supset \cdots
\end{equation*}
(where $X_i=t^i(X)$) and proved the main theorem that the tight span $T(X)$ of any metric space $X=(X,d)$ can be expressed as
    \begin{equation*}
      T(X)=\uptau \cup\overline{C(X)} \quad \text{with} \quad \uptau \cap \overline{C(X)}=\overline{C(X)}_{*},
    \end{equation*}
where $\uptau=\displaystyle\cap_{i\geq1}T(X_i)$ and $\overline{C(X)}_{*}$ is a certain subset of $\overline{C(X)}$ called the roots. (For details see \cite{Turaev2018}; we use a special case of this theorem where $X$ is a finite metric space.)

In Section 2, we characterize trim subspaces of the Manhattan plane ($\mathbb{R}^2_1$), define the metric centers of a triple of points and of a finite subspace $X$ of $\mathbb{R}^2_1$ and identify the abstract space $t(X)$ with the metric center set $m(X)$ which lives again in $\mathbb{R}^2_1$.

In Section 3, we give a simple method to obtain the metric center set $m(X)$ and in Section 4 we explain how the (metric quotient of the) trimming cylinder can be embedded into $\mathbb{R}^2_1$.

In Section 5, we describe a conceptually lucid way of constructing the tight span of a finite subset of $\mathbb{R}^2_1$ as an application of the related theorem of Turaev and we give an algorithm to implement it. (For another algorithm to construct the tight span of a finite subset of $\mathbb{R}^2_1$ see \cite{Kilic2021}.)

\section{Trimming of Finite Subspaces of the Manhattan Plane}

We remind that the Manhattan plane is the metric space $(\mathbb{R}^2,d_1)$ with \linebreak $d_1((x_1,y_1),(x_2,y_2))=|x_1-x_2|+|y_1-y_2|$ for $(x_1,y_1), (x_2,y_2)\in \mathbb{R}^2$. We denote this space by $\mathbb{R}^2_1$.

We first note the following simple and useful property for the relation of \linebreak ``betweenness" in the Manhattan plane.
\begin{lemma}\label{lem1}
  A point $(x_0,y_0)\in \mathbb{R}^2_1$ lies between $(x_1,y_1)$ and $(x_2,y_2)\in \mathbb{R}^2_1$ in the sense of Menger if and only if $x_0$ lies between $x_1$ and $x_2$ (in $\mathbb{R}$) and $y_0$ lies between $y_1$ and $y_2$ ($x_1\leq x_0 \leq x_2$ or $x_2 \leq x_0 \leq x_1$, respectively for $y_0$).
\end{lemma}
Recall that $\triangle_{(x_0,y_0),(x_1,y_1)(x_2,y_2)}=0$ if and only if $(x_0,y_0)$ lies between $(x_1,y_1)$ and $(x_2,y_2)$.

Now, we give a criterion for a finite subset of $\mathbb{R}^2_1$ to be trim:
\begin{proposition}\label{theo0}
   Let $(X,d_1)\subset\mathbb{R}^2_1$ be a finite subspace with $|X|\geq4$. Let $R_X$ denote the minimal rectangle in $\mathbb{R}^2_1$ containing $X$ with edges parallel to the axes. Then,  $X$ is trim if and only if each edge of $R_X$ contains at least two points of $X$.
\end{proposition}
\begin{proof}
 First assume that $X$ is trim. Without loss of generality consider the right edge of $R_X$. There must be at least one point of $X$ on this edge, since otherwise the rectangle would not be minimal. Assume, there is only one point of $X$ on this edge, say $(x_0,y_0)$. If  $(x_1,y_1)$ and $ (x_2,y_2)$ are two other points of $X$, then necessarily $x_1 < x_0$ and $x_2 < x_0$ so that by Lemma \ref{lem1}, $\triangle_{(x_0,y_0),(x_1,y_1)(x_2,y_2)}>0$ and then
  $\underline{d}_1((x_0,y_0))>0$ since $X$ is finite. This contradicts the assumption that $X$ is trim. So, there must be at least two points of $X$ on the considered edge of $R_X$.

  Conversely, let us assume that each edge of $R_X$ contains at least two points of $X$. Now, a point of $X$ can be on the boundary or inside of $R_X$. First, consider a point $(x_0,y_0)\in X$ lying on an edge of $R_X$, say without loss of generality, on the right edge. Then choose another point $(x_0,y_1)\in X$ on the same edge and a third point $(x_2,y_2)\in X$ on the top or bottom edge of $R_X$ such that $y_0$ lies between $y_1$ and $y_2$ (see Figure \ref{fig3}). Then by the Lemma \ref{lem1}, $\triangle_{(x_0,y_0), (x_0,y_1)(x_2,y_2)}=0$ and then $\underline{d}_1((x_0,y_0))=0$.

     Now, let us assume that $(x_0,y_0)\in X$ lies in the interior of $R_X$. Divide the rectangle $R_X$ with four subrectangles by the vertical and horizontal lines passing through $(x_0,y_0)$ (see Figure \ref{fig12}). There must be ``boundary" points of $X$ (i.e. points of $X$ lying on the boundary of $R_X$) either on both of the subrectangles $\rm{I}$ and $\rm{III}$, or on both of the subrectangles $\rm{II}$ and $\rm{IV}$. Because otherwise there would be neighbouring subrectangles (such as $\rm{I}$ and $\rm{II}$, or $\rm{II}$ and $\rm{III}$, or $\rm{III}$ and $\rm{IV}$, or $\rm{IV}$ and $\rm{I}$) containing no boundary points, which would contradict the minimality of $R_X$. Thus by Lemma \ref{lem1}, we have $\underline{d}_1((x_0,y_0))=0$ and $X$ is trim.
  \end{proof}
\begin{figure}
  \centering
 \scalebox{1}{ \begin{tikzpicture}
    \draw (0,0) -- (5,0);     \draw (0,0) -- (0,3);      \draw (0,3) -- (5,3);      \draw (5,3) -- (5,0);
    \draw [fill] (0,0) circle (.05); \draw [fill] (0,1) circle (.05);  \draw [fill] (0,1.5) circle (.05);  \draw [fill] (1,3) circle (.05);
    \draw [fill] (3.5,3) circle (.05); \draw [fill] (5,1.5) circle (.05);  \draw [fill] (5,.7) circle (.05);  \draw [fill] (4,0) circle (.05);
    \draw [fill] (3,2) circle (.05); \draw [fill] (3.5,.5) circle (.05);  \draw [fill] (1,.3) circle (.05);
    \node [above] at (3.5,3) {$(x_2,y_2)$};     \node [right] at (5,.7) {$(x_0,y_1)$};       \node [right] at (5,1.5) {$(x_0,y_0)$};
\end{tikzpicture}}
\caption{An example of the positions of $(x_0,y_0)$, $(x_0,y_1)$ and $(x_2,y_2)$.}\label{fig3}
  \end{figure}
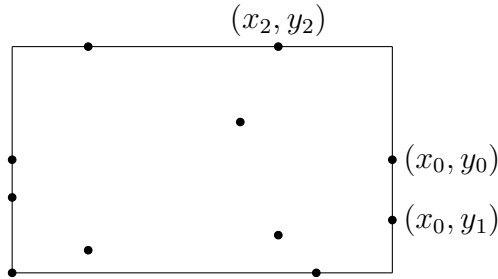
  \begin{figure}
  \centering
  \begin{tikzpicture} %bölgeler
    \draw (0,0) -- (5,0);  \draw (5,3) -- (5,0);  \draw (0,3) -- (5,3);  \draw (0,0) -- (0,3);
     \draw [fill] (3,2) circle (.05);  \draw (-0.2,2) -- (5.2,2);  \draw (3,-0.2) -- (3,3.2);
     \node [below left] at (3,2) {$(x_0,y_0)$};
     \node [above] at (3.5,2.1) {$\rm{I}$}; \node [above] at (2,2.1) {$\rm{II}$};  \node [below] at (2,1) {$\rm{III}$};  \node [below] at (3.5,1) {$\rm{IV}$};
\end{tikzpicture}
\caption{The four subrectangles of $R_X$.}\label{fig12}
  \end{figure}
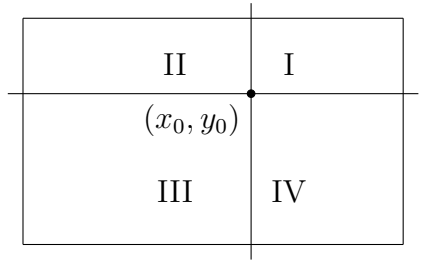
\begin{remark}\label{remark1}
  The last part of the above proof shows that points $(x_0,y_0)\in X$ lying in the interior of $R_X$ are always in-between points (i.e. $\underline{d}_1((x_0,y_0))=0$), since one point on each edge of $R_X$ is sufficient for this argument and this is obviously fulfilled by minimality of $R_X$.
\end{remark}

Our next goal is to give an interpretation of the process of trimming of finite subsets of $\mathbb{R}^2_1$ in terms of metric centers of this subset. We first give some auxiliary definitions.

For a two point subset $\{a,b\}\subset \mathbb{R}^2_1$ let us denote $R_{\{a,b\}}$ by $R_{ab}(=R_{ba})$, which is the minimal rectangle containing $a$ and $b$, with edges parallel to the axes. In view of Lemma \ref{lem1}, a point $c\in \mathbb{R}^2_1$ is between $a$ and $b$ if and only if $c\in R_{ab}$.

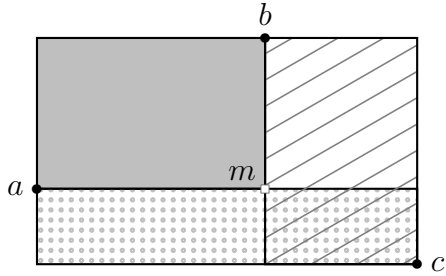
\begin{figure}
\centering
\begin{pspicture}(-1,-1)(6,4)
%\pspolygon(0,0)(5,0)(5,3)(0,3)(0,0)
%\pspolygon[fillstyle=vlines,hatchsep=1pt,fillcolor=lightgray](3,1)(5,1)(5,3)(3,3)
\pspolygon[fillstyle=solid,hatchwidth=0.7pt,fillcolor=lightgray](0,1)(3,1)(3,3)(0,3)
\pspolygon[fillstyle=dots,hatchwidth=0.6pt,dotsize=1pt,hatchsep=3pt,hatchcolor=lightgray](0,0)(5,0)(5,1)(0,1)
\pspolygon[fillstyle=hlines,hatchwidth=0.6pt,hatchsep=8pt,hatchangle=30,hatchcolor=gray](3,0)(5,0)(5,3)(3,3)
\uput[l](0,1){$a$} \psdot(0,1)
\uput[u](3,3){$b$} \psdot(3,3)
\uput[r](5,0){$c$} \psdot(5,0)
\uput[ul](3,1){$m$} \psdot[dotstyle=square](3,1) %,fillcolor=black
\end{pspicture}
\caption{Intersection of the rectangles $R_{ab}$, $R_{ac}$ and $R_{bc}$.}\label{fig4}
\end{figure}

For three points $a, b, c \in \mathbb{R}^2_1$, the intersection of three rectangles $R_{ab}$, $R_{ac}$ and $R_{bc}$ is a singleton (see Figure \ref{fig4}) and it is called the metric center of $\{a, b, c\}$, or in other words:
\begin{definition}
  Let $a,b,c \in \mathbb{R}^2_1$. The unique point $m\in \mathbb{R}^2_1$ satisfying
  \begin{align*}
d_1(a,b)&=d_1(a,m)+d_1(m,b),\\
d_1(a,c)&=d_1(a,m)+d_1(m,c),\\
d_1(b,c)&=d_1(b,m)+d_1(m,c)
\end{align*}
is called the metric center of $\{a, b, c\}$. (We denote it sometimes by $m=m(a,b,c)$.)
\end{definition}

As a consequence of these equations one gets $d_1(a,m)=\triangle_{a,bc}$, $d_1(b,m)=\triangle_{b,ac}$ and $d_1(c,m)=\triangle_{c,ab}$.
%%% 그래픽 추가 fig4

\begin{definition}\label{def1}
 Let $X\subset \mathbb{R}^2_1$ be a finite subspace and $a\in X$. Assume that $\triangle_{a,bc}$ (with $b,c\in X$) is the minimal Gromov product at $a$. Then, the metric center of $\{a,b,c\}$ is called the metric center associated with $a$ (or simply, of $a$) and is denoted by $m_a$ (or sometimes by $m(a)$). It holds then $\underline{d}_1(a)=\min_{\substack{u,v\in X\setminus\{a\}\\ u\neq v}} \triangle_{a,uv}=\triangle_{a,bc}=d_1(a,m_a)$. The set of all metric centers of the elements of $X$ is called the metric center of $X$ and it is denoted by $m(X)$.

 For reasons of convenience, we define the metric center $m(X)$ of a two-point space $X=\{a,b\}\subset (\mathbb{R}^2,d_1)$ as $m(X)=\{\frac{a+b}{2}\}$ and the metric center $m(X)$ of a singleton $X=\{a\}\subset (\mathbb{R}^2,d_1)$ as $m(X)=X$.
\end{definition}

\begin{remark}\label{remark2}
  If the minimal Gromov product at $a$ is realized simultaneously by another triple $\{a, d, e\}$; then it can be easily seen that $m(a,b,c)=m(a,d,e)$ so that $m_a$ is well-defined. The following useful properties hold:
  \begin{equation*}
    d_1(a,x)=d_1(a,m_a)+d_1(m_a,x) \;\mbox{for any}\; x\in X\setminus\{a\}
  \end{equation*}
  and
    \begin{equation*}
    d_1(a,m_x)=d_1(a,m_a)+d_1(m_a,m_x) \;\mbox{for any}\; x\in X.
  \end{equation*}
(These properties will be clear by the simple geometric construction of the metric center $m(X)$ we give in the next section.)
\end{remark}

 Now consider $x,y \in X \subset \mathbb{R}^2_1$ and assume that $\overline{x}=\overline{y}\in t(X)$. By definition,
  \[d_1 (x,y)=\underline{d}_1(x)+\underline{d}_1(y)\]
  and thus
  \[d_1 (x,y)= d_1(x,m_{x})+d_1(y,m_{y}).\]
  On the other hand,
  \begin{align*}
  d_1(x,y) & = d_1(x,m_{x})+ d_1(m_{x},y) \\
   &= d_1 (x,m_x)+d_1(y,m_y)+d_1(m_y,m_x)
\end{align*}
by the above properties, so that we get $d_1(m_y,m_x)=0$ and thus $m_x=m_y$. In other words, if two elements are identified during the trimming, then they have the same metric centers. This interesting property enables us to define a map $t(X)\to m(X)$, which turns out to be an isometry.

\begin{theorem}\label{theo1}
  Let $X\subset(\mathbb{R}^2,d_1)$ be a finite subspace and let $m_x$ be the metric center of $x\in X$. Then the map
  \begin{align*}
   f: (t(X), d^1) &\rightarrow (m(X),d_1) \\
    \overline{x} &\mapsto  m_x
  \end{align*}
  is an isometry.
\end{theorem}
\begin{proof}
  $f$ is obviously surjective. By the definition of $d^1$ and the above properties of centers we can write
\begin{align*}
  d^1(\overline{x},\overline{y})&=d_1(x,y)- \underline{d}_1(x)- \underline{d}_1(y)\\
  &= d_1(x,y)- d_1(x,m_{x})- d_1(y,m_{y}) \\
  &= d_1(x,m_{x})+d_1(m_{x},y)- d_1(x,m_{x})- d_1(y,m_{y})\\
  &=d_1(m_x,y)-d_1(y, m_y)\\
  &= d_1(m_{x},m_{y})+d_1(m_{y},y)- d_1(y,m_{y}) \\
  &= d_1(m_{x},m_{y}),
\end{align*}
so that $f$ is an isometry.
 \end{proof}

  This relationship can also be expressed as a commutative diagram

   \begin{figure}[h]
   \centering
 \includegraphics[scale=1]{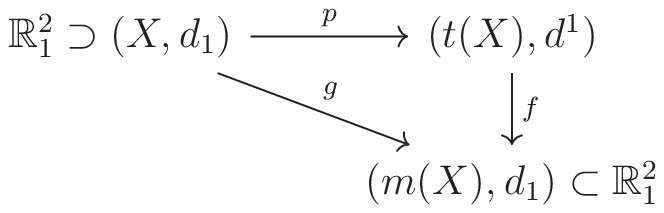}
 \end{figure}
\hspace{-0.8cm} where $p$ is the trimming projection and $g$ is the operation of taking the centers (or, the ``center projection"). This means that the abstract trimming story of a finite subset of the Manhattan plane can be staged in this plane itself.

Instead of the trimming sequence
\[
(X,d)\xrightarrow{p}(t(X),d^1)\xrightarrow{p_1}(t^2(X),d^2)\xrightarrow{p_2} \cdots
\]
we can work with the embedded sequence
\[
(X,d)\xrightarrow{g} (m(X),d_1)\xrightarrow{g_1} (m(m(X)),d_1)\xrightarrow{g_2}\cdots ,
\]
where all terms live in $\mathbb{R}^2_1 $. We will denote the $i^{th}$ iterate $m(m\cdots(m(X)))$ by $m^i(X)$, the center projection $(m^i(X),d_1) \to(m^{i+1}(X),d_1)$ by $g_i$ and call the sequence
\begin{equation*}
\resizebox{.97\hsize}{!}{$(X,d)\xrightarrow{g} (m(X),d_1)\xrightarrow{g_1} (m^2(X),d_1)\to\cdots \to (m^i(X),d_1) \xrightarrow{g_i}(m^{i+1}(X),d_1)\to \cdots $}
\end{equation*}
the metric center sequence of $(X,d)\subset \mathbb{R}^2_1$.

%%%%%%.
\section{How to Find the Metric Center of a Finite Subspace of the Manhattan Plane?}

We now want to give a simple geometric construction to obtain the metric center $m(X)\subset \mathbb{R}_1^2$ of a finite subspace $(X,d_1)\subset \mathbb{R}_1^2$. Let $X$ have $n$ points and let us order the abcissas of these points as
\[x_1\leq x_2\leq \cdots\leq x_{n-1}\leq x_n,\]
and the ordinates as
\[y_1\leq y_2 \leq \cdots \leq y_{n-1}\leq y_n.\]

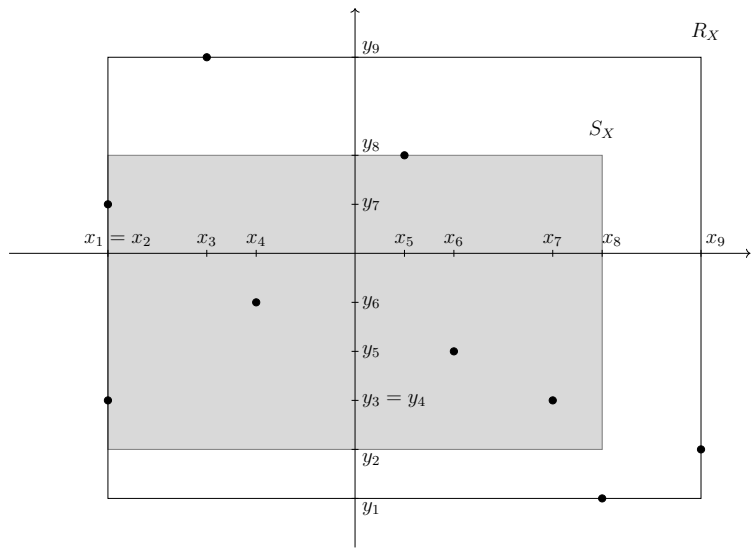
\begin{figure}[h]
\centering
\scalebox{.65}{
\begin{tikzpicture} %[x={10.0pt},y={10.0pt}]

    \draw [fill=gray!30, draw=gray] (-5,2)--(5,2)--(5,-4)--(-5,-4)--cycle;
      \draw [draw=black] (-5,4)--(7,4)--(7,-5)--(-5,-5)--cycle;
         \foreach \x in {-5,-3,-2,1,2,4,5,7} {%
    \draw ($(\x,0) + (0,-\TickSize)$) -- ($(\x,0) + (0,\TickSize)$)
        node [above] {};
}

\foreach \y in {-5,-4,-3,-2,-1,1,2,4} {%
    \draw ($(0,\y) + (-\TickSize,0)$) -- ($(0,\y) + (\TickSize,0)$)
        node [right] {};
}
%\node[fill=black,star,star points=5,star point ratio=0.5] at (-3,2) {};   \node[fill=black,star,star points=5,star point ratio=0.5] at (5,-4) {};  % \node [below] at (5,-4.2) {$m_b=m_c$};

\node [above] at (7.1,4.2) {$R_X$};  \node [above] at (5,2.2) {$S_X$};
 \draw[black,fill=black](-5,1) circle (.4ex); \draw[black,fill=black](-3,4) circle (.4ex); \draw[black,fill=black](1,2) circle (.4ex);\draw[black,fill=black](-5,-3) circle (.4ex);
  \draw[black,fill=black](-2,-1) circle (.4ex); \draw[black,fill=black](2,-2) circle (.4ex); \draw[black,fill=black](4,-3) circle (.4ex);
 \draw[black,fill=black](5,-5) circle (.4ex);\draw[black,fill=black](7,-4) circle (.4ex);

      \draw [->] (-7,0) -- (8,0);  \draw [->] (0,-6) -- (0,5);
     \node [above] at (-4.8,0) {$x_1=x_2$};  \node [above] at (-3,0) {$x_3$}; \node [above] at (-2,0) {$x_4$};  \node [above] at (1,0) {$x_5$};
     \node [above] at (2,0) {$x_6$}; \node [above] at (4,0) {$x_7$};  \node [above] at (5.2,0) {$x_8$}; \node [above] at (7.3,0) {$x_9$};
    % \node [above] at (-3,4) {$a$}; \node [below] at (5,-5) {$b$}; \node [below right] at (7,-4) {$c$};

     \node [right] at (0,-5.2) {$y_1$}; \node [right] at (0,-4.2) {$y_2$}; \node [right] at (0,-3) {$y_3=y_4$}; \node [right] at (0,-2) {$y_5$}; \node [right] at (0,-1) {$y_6$};
     \node [right] at (0, 1) {$y_7$}; \node [right] at (0,2.2) {$y_8$}; \node [right] at (0,4.2) {$y_9$};
\end{tikzpicture}}
\caption{An example of the minimal rectangle $R_X$ and the secondary rectangle $S_X$ for a 9-point subspace $X\subset \mathbb{R}^2_1$ with $X=\{(-5,-3),(-5,1),(-3,4),$ $(-2,-1),(1,2),(2,-2),(4,-3),(5,-5),(7,-4)\}$.}\label{fignew0}
\end{figure}

We then define a secondary rectangle $S_X=[x_2,x_{n-1}]\times [y_2, y_{n-1}]$ (see Figure \ref{fignew0}). Recall that $R_X$ was the minimal, axes-parallel rectangle containing $X$, which can be expressed with these notations as $R_X=[x_1,x_n]\times [y_1,y_n]$.

The following proposition gives a nice device to determine the metric center $m(X):$
\begin{proposition}
  The metric center $m_a$ of a point $a\in X$ is the point of $S_X$ nearest to $a$ (with respect to the $d_1-$metric, or what the same is, with respect to the Euclidean metric; we call this point the projection of $a$ on $S_X$). Especially, for a point $a\in X$ belonging to $S_X$, we get $m_a=a$. The metric center $m(X)$ of $X$ thus consists of the points of $X$ contained in $S_X$ and the projections of the points of $X$ outside $S_X$ onto $S_X$.
\end{proposition}
\begin{proof}
In Figure \ref{fignew101} all possible configurations of $R_X$ (the minimal axes-parallel rectangle containing $X$), the secondary rectangle $S_X$ and the points of $X$ lying outside $S_X$ are depicted. In the first row the possible relative positions of $R_X$ and $S_X$ (up to ``symmetry") are shown and in the corresponding columns below the possible placements of points of $X$ lying outside $S_X$ (again up to ``symmetry") are shown.

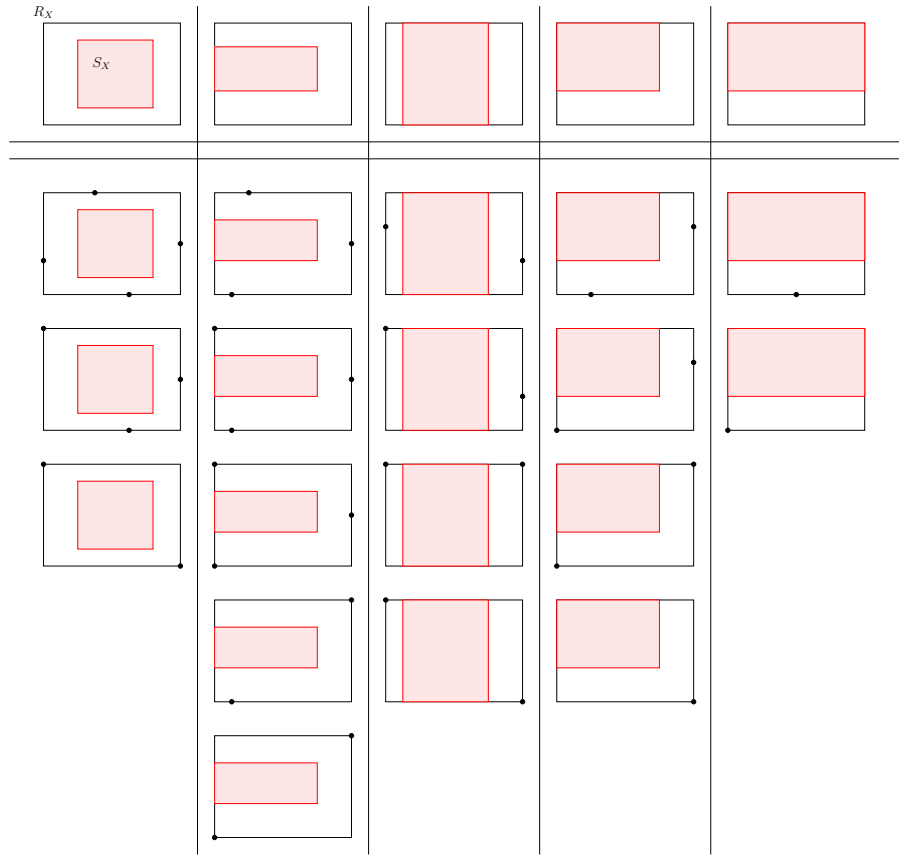
\begin{figure}[h]
  \centering
   \scalebox{.45}{\begin{tikzpicture}
   %%first colon and lines % \node[above] at (3,4) {$a$};
\draw[black] (-1,-0.5) -- (25,-0.5); \draw[black] (-1,-1) -- (25,-1); \draw[black] (4.5,3.5) -- (4.5,-21.5);
\node[above] at (0,3) {$R_X$};
  \draw [black](0,0) rectangle (4,3); \draw [red,fill=red!10](1,0.5) rectangle (3.2,2.5);
  \draw [black](0,-5) rectangle (4,-2); \draw [red,fill=red!10](1,-4.5) rectangle (3.2,-2.5); \draw[black,fill=black](1.5,-2) circle (.35ex); \draw[black,fill=black](4,-3.5) circle (.35ex); \draw[black,fill=black](2.5,-5) circle (.35ex); \draw[black,fill=black](0,-4) circle (.35ex);
  \draw [black](0,-9) rectangle (4,-6); \draw [red,fill=red!10](1,-8.5) rectangle (3.2,-6.5); \draw[black,fill=black](0,-6) circle (.35ex); \draw[black,fill=black](4,-7.5) circle (.35ex); \draw[black,fill=black](2.5,-9) circle (.35ex);
  \draw [black](0,-13) rectangle (4,-10); \draw [red,fill=red!10](1,-12.5) rectangle (3.2,-10.5); \draw[black,fill=black](0,-10) circle (.35ex); \draw[black,fill=black](4,-13) circle (.35ex);
 \node[above] at (1.7,1.5) {$S_X$};
    %%second colon and lines
 \draw[black] (9.5,3.5) -- (9.5,-21.5);
  \draw [black](5,0) rectangle (9,3); \draw [red,fill=red!10](5,1) rectangle (8,2.3);
  \draw [black](5,-5) rectangle (9,-2); \draw [red,fill=red!10](5,-4) rectangle (8,-2.8); \draw[black,fill=black](6,-2) circle (.35ex); \draw[black,fill=black](9,-3.5) circle (.35ex); \draw[black,fill=black](5.5,-5) circle (.35ex);
  \draw [black](5,-9) rectangle (9,-6); \draw [red,fill=red!10](5,-8) rectangle (8,-6.8); \draw[black,fill=black](5,-6) circle (.35ex); \draw[black,fill=black](9,-7.5) circle (.35ex); \draw[black,fill=black](5.5,-9) circle (.35ex);
  \draw [black](5,-13) rectangle (9,-10); \draw [red,fill=red!10](5,-12) rectangle (8,-10.8); \draw[black,fill=black](5,-10) circle (.35ex); \draw[black,fill=black](9,-11.5) circle (.35ex); \draw[black,fill=black](5,-13) circle (.35ex);
  \draw [black](5,-17) rectangle (9,-14); \draw [red,fill=red!10](5,-16) rectangle (8,-14.8); \draw[black,fill=black](9,-14) circle (.35ex); \draw[black,fill=black](5.5,-17) circle (.35ex);
  \draw [black](5,-21) rectangle (9,-18); \draw [red,fill=red!10](5,-20) rectangle (8,-18.8);  \draw[black,fill=black](9,-18) circle (.35ex); \draw[black,fill=black](5,-21) circle (.35ex);

     %%third colon and lines
 \draw[black] (14.5,3.5) -- (14.5,-21.5);
  \draw [black](10,0) rectangle (14,3); \draw [red,fill=red!10](10.5,0) rectangle (13,3);
  \draw [black](10,-5) rectangle (14,-2); \draw [red,fill=red!10](10.5,-5) rectangle (13,-2); \draw[black,fill=black](10,-3) circle (.35ex); \draw[black,fill=black](14,-4) circle (.35ex);
  \draw [black](10,-9) rectangle (14,-6); \draw [red,fill=red!10](10.5,-9) rectangle (13,-6); \draw[black,fill=black](10,-6) circle (.35ex); \draw[black,fill=black](14,-8) circle (.35ex);
  \draw [black](10,-13) rectangle (14,-10); \draw [red,fill=red!10](10.5,-13) rectangle (13,-10); \draw[black,fill=black](10,-10) circle (.35ex); \draw[black,fill=black](14,-10) circle (.35ex);
  \draw [black](10,-17) rectangle (14,-14); \draw [red,fill=red!10](10.5,-17) rectangle (13,-14); \draw[black,fill=black](10,-14) circle (.35ex); \draw[black,fill=black](14,-17) circle (.35ex);

 %%fourth colon and lines
 \draw[black] (19.5,3.5) -- (19.5,-21.5);
  \draw [black](15,0) rectangle (19,3); \draw [red,fill=red!10](15,1) rectangle (18,3);
  \draw [black](15,-5) rectangle (19,-2); \draw [red,fill=red!10](15,-4) rectangle (18,-2); \draw[black,fill=black](16,-5) circle (.35ex); \draw[black,fill=black](19,-3) circle (.35ex);
  \draw [black](15,-9) rectangle (19,-6); \draw [red,fill=red!10](15,-8) rectangle (18,-6); \draw[black,fill=black](15,-9) circle (.35ex); \draw[black,fill=black](19,-7) circle (.35ex);
  \draw [black](15,-13) rectangle (19,-10); \draw [red,fill=red!10](15,-12) rectangle (18,-10); \draw[black,fill=black](15,-13) circle (.35ex); \draw[black,fill=black](19,-10) circle (.35ex);
  \draw [black](15,-17) rectangle (19,-14); \draw [red,fill=red!10](15,-16) rectangle (18,-14); \draw[black,fill=black](19,-17) circle (.35ex);

   %%fifth colon and lines
  \draw [black](20,0) rectangle (24,3); \draw [red,fill=red!10](20,1) rectangle (24,3);
  \draw [black](20,-5) rectangle (24,-2); \draw [red,fill=red!10](20,-4) rectangle (24,-2); \draw[black,fill=black](22,-5) circle (.35ex);
  \draw [black](20,-9) rectangle (24,-6); \draw [red,fill=red!10](20,-8) rectangle (24,-6); \draw[black,fill=black](20,-9) circle (.35ex);

\end{tikzpicture}}
\caption{Possible configurations of $R_X, S_X$ and the points of $X$ lying outside $S_X$.}\label{fignew101}
  \end{figure}

We first remark that by the minimality of $R_X$ there exists at least one point on every edge of $R_X$ belonging to $X$ and all points of $X$ lying in the interior of $R_X$ are in-between points (i.e. have zero pendant length), so that their metric centers coincide with them.

We next remark that, in cases one edge of $S_X$ is contained in an edge of $R_X$, there must be at least two points of $X$ on this edge of $R_X$ by the definition of $S_X$ and in such a case, points of $X$ lying in the corresponding edge of $S_X$ must be in-between points since one can use auxiliary points on approximate neighbouring edges of $R_X$ to see this. Thus, these points also coincide with their metric centers.

Now, there remain the points of $X$, which lie on the boundary of $R_X$ but do not belong to $S_X$. For these points there are two typical cases as shown in Figure \ref{fignew100}, whereby the rays in Figure \ref{fignew100}(a) and the line in Figure \ref{fignew100}(b) are determined by the relevant edges of $S_X$.

In Figure \ref{fignew100}(a), the metric center $m_a$ of the corner $a$ is the corner where the rays emanate: Either there are points $b$, $c$ on the rays belonging to $X$, in which case $m_a$ is the metric center of the triple $(a, b, c)$; or the origin of the rays belongs to $X$, in which case one can take any third auxiliary point of $X$ and again the origin of the rays becomes the metric center $m_a$. (If one chooses any other two points of $X$, then their center is either farther away from $a$, or the same corner again.)

  \begin{figure}[h]
\centering
    \scalebox{.65}{\begin{tikzpicture}%[x={10.0pt},y={10.0pt}]
   %firstrowrectangles
   %1st
     \draw [fill=red!10](1,2) rectangle (4,0);
    \draw[red,->] (1,2) -- (1,-1);  \draw[red,dotted] (1,2) -- (0,3);  \draw[red,fill=red](1,2) circle (.3ex); \node[above] at (1,2) {$m_a$};
    \draw[red,->] (1,2) -- (5,2);
    \draw[black,fill=black](0,3) circle (.3ex); \node[above] at (0,3) {$a$};
  \draw [black](0,0) rectangle (4,3);
%2nd
\draw [fill=red!10](7,0) rectangle (10,2); \draw [blue,dotted](6,3) rectangle (7,1);  \draw [blue,dotted](6,3) rectangle (9,2);  \draw [blue,dotted](9,2) rectangle (7,1);
\draw[red,->] (7,2) -- (7,-1); \draw[red,->] (7,2) -- (11,2); \draw[red,fill=red](7,2) circle (.3ex); \node[above] at (7,2) {$m(a,b,c)$};
 \draw[black,fill=black](6,3) circle (.3ex); \node[above] at (6,3) {$a$};
  \draw[black,fill=black](7,1) circle (.3ex); \node[below left] at (7,1) {$b$};
   \draw[black,fill=black](9,2) circle (.3ex); \node[above right] at (9,2) {$c$};
   \draw [black](6,0) rectangle (10,3);
%3rd
\draw [fill=red!10](13,0) rectangle (16,2); \draw [blue,dotted](12,3) rectangle (14,0.5);  \draw [blue,dotted](12,3) rectangle (15.5,1);  \draw [blue,dotted](15.5,1) rectangle (14,0.5);
\draw[red,->] (13,2) -- (13,-1); \draw[red,->] (13,2) -- (17,2); \draw[red,fill=red](14,1) circle (.3ex); \node[above] at (14,1) {$m(a,u,v)$};
 \draw[black,fill=black](12,3) circle (.3ex); \node[above] at (12,3) {$a$};
  \draw[black,fill=black](14,0.5) circle (.3ex); \node[below left] at (14,0.5) {$u$};
   \draw[black,fill=black](15.5,1) circle (.3ex); \node[above right] at (15.5,1) {$v$};
    \draw [black](12,0) rectangle (16,3);
%% (a) ve (b)
\node[] at (7.5,-2) {(a)};  \node[] at (7.5,-8) {(b)};
   %secondrowrectangles
   %1st
   \draw [fill=red!10](0,-7) rectangle (4,-5);
    \draw[red,<->] (-1,-5) -- (4.5,-5);  \draw[red,dotted] (2,-4) -- (2,-5);  \draw[red,fill=red](2,-5) circle (.3ex);  \node[above left] at (2,-5) {$m_a$};
     \draw[black,fill=black](2,-4) circle (.3ex); \node[above] at (2,-4) {$a$};
  \draw [black](0,-7) rectangle (4,-4);
%2nd
 \draw [fill=red!10](6,-7) rectangle (10,-5); \draw [blue,dotted](7,-5) rectangle (8,-4);  \draw [blue,dotted](7,-5) rectangle (9,-6.5);  \draw [blue,dotted](8,-4) rectangle (9,-6.5);
 \draw[red,<->] (5.5,-5) -- (10.5,-5); \draw[red,fill=red](8,-5) circle (.3ex);  \node[below right] at (8,-5) {$m(a,b,c)$};
 \draw[black,fill=black](8,-4) circle (.3ex); \node[above] at (8,-4) {$a$};
 \draw[black,fill=black](7,-5) circle (.3ex); \node[below left] at (7,-5) {$b$};
 \draw[black,fill=black](9,-6.5) circle (.3ex); \node[right] at (9,-6.5) {$c$};
 \draw [black](6,-7) rectangle (10,-4);
%3rd
 \draw [fill=red!10](12,-7) rectangle (16,-5); \draw [blue,dotted](12.5,-6) rectangle (14,-4);  \draw [blue,dotted](12.5,-6) rectangle (15,-6.5);  \draw [blue,dotted](14,-4) rectangle (15,-6.5);
\draw[red,<->] (11.5,-5) -- (16.5,-5); \draw[red,fill=red](14,-6) circle (.3ex);  \node[above] at (14,-6) {$m(a,u,v)$};
 \draw[black,fill=black](14,-4) circle (.3ex); \node[above] at (14,-4) {$a$};
 \draw[black,fill=black](12.5,-6) circle (.3ex); \node[left] at (12.5,-6) {$u$};
 \draw[black,fill=black](15,-6.5) circle (.3ex); \node[right] at (15,-6.5) {$v$};
 \draw [black](12,-7) rectangle (16,-4);
\end{tikzpicture}
}

\caption{Two typical cases of points $a\in X$ lying on the boundary of $R_X$ which do not belong to $S_X$.} \label{fignew100}
\end{figure}
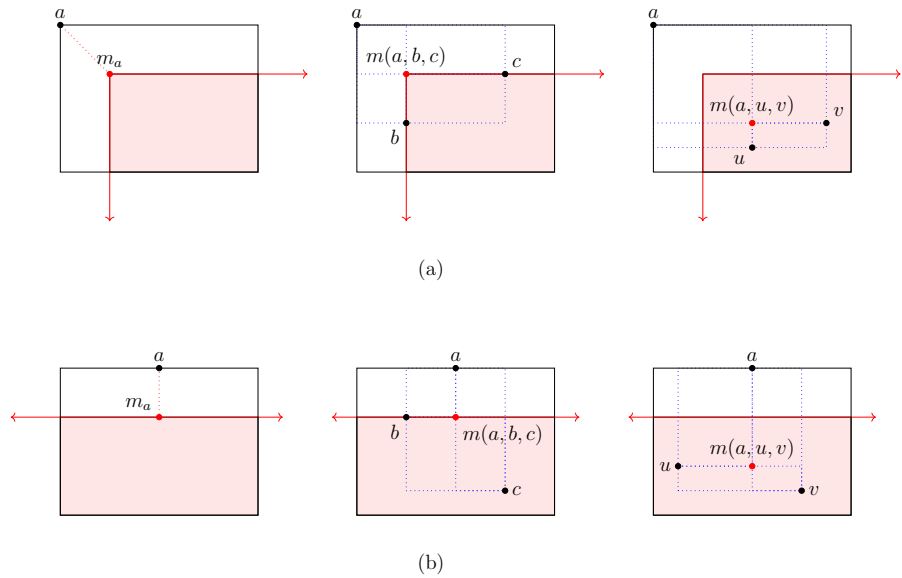

In Figure \ref{fignew100}(b), there is at least one point of $X$ lying on the line and choosing an auxiliary point of $X$ lying on the other part of the vertical line passing through $a$ one obtains $m_a$ which is the projection of $a$ onto the line. (If one chooses any other two points of $X$, then their center is either farther away from $a$, or the same projection point again.)
\end{proof}

%%%%%%%%
\section{How to Embed the Trimming Cylinder into the Manhattan Plane?}
We have seen in Theorem \ref{theo1} that the space $(t(X),d^1)$, obtained by trimming a finite subspace $(X,d_1)\subset \mathbb{R}^2_1$, can be isometrically embedded into the $\mathbb{R}^2_1$, by the canonical map $\overline{x}\mapsto m_x$, sending a class to the metric center of a representant of it. This enables us to work with the embedded version of the abstract trimming sequence, which we called the metric center sequence. This procedure can be extended to the trimming cylinder also. We can define the corresponding metric center cylinder $MC(X)$ as a pseudo-metric graph with the vertex set $\bigsqcup_{i\geq 0}m^i(X)$ and edges between $a_ i\in m^i(X)$ and $g_i(a_i)\in m^{i+1}(X)$ with weight $\underline{d}_1(a_i)$. It is instructive to look at the first portion $MC_1(X)$ of this cylinder (which repeats itself in the following stages since in every stage the starting space is a subspace of $\mathbb{R}^2_1$):

  \begin{figure}[h]
   \centering
   \begin{tikzcd}
    & (t(X),d^1) \arrow{d}{f}\\
\mathbb{R}^2_1 \supset(X,d_1) \arrow{r}{g} \arrow{ru}{p} & (m(X),d_1)\subset \mathbb{R}^2_1
\end{tikzcd}
 \end{figure}

The edge between $x\in X$ and $p(x)\in t(X)$ in $C_1(X)$ has weight $\underline{d}_1(x)$; the edge between $x\in X$ and $g(x)\in m(X)$ has also weight $\underline{d}_1(x)$. On the other hand, the distance $d_1(x,g(x))=d_1(x, m_x)$ also equals $\underline{d}_1(x)$ so that one can embed this ``abstract" edge as the straight line segment connecting $x$ with $m_x$ in $\mathbb{R}^2_1$. In case $\underline{d}_1(x)=0$ we have $m_x=x$ so that the metric quotient $\overline{MC_1(X)}$ consists exactly of the set $m(X)$ together with the segments $[x,m_x]$ in $\mathbb{R}^2_1$. We note that the segments $[x,m_x]$ and $[y,m_y]$ in $\mathbb{R}^2_1$ cannot cross themselves (unless maybe in the special case $m_x=m_y$, where they are concatenated) since by Remark \ref{remark2} \[d_1(x,y)=d_1(x,m_x)+d_1(m_x,y)=d_1(x,m_x)+d_1(m_x,m_y)+d_1(m_y,y)\] holds, which means that the points $x,m_x,m_y$ and $y$ lie on a geodesic (shortest path) between $x$ and $y$. (Likewise in the following step $x,m_x,m^2_x,m^2_y,m_y$ and $y$ will lie on a geodesic etc.) Iterating this procedure, which stabilizes after finitely many steps, one obtains an embedding of $\overline{MC(X)}$ into $\mathbb{R}^2_1$ for the finite subspace $X\subset \mathbb{R}^2_1$. We illustrate this process by the following:

\begin{example}\label{example1}
Let $X$ be as in Figure \ref{fignew0}, i.e. $X=\{ a=(-5,-3),b=(-5,1),c=(-3,4),d=(-2,-1),e=(1,2),$ $f=(2,-2),g=(4,-3),h=(5,-5), i=(7,-4)\}$.

We show in Figure \ref{figex1} the successive metric centers $m(X), m^2(X)$ and $m^3(X)$ in the left column of the figure.
$m^3(X)$ is a trim space so that the metric center sequence stabilizes at this stage (the minimal rectangle of $m^3(X)$ is shown by the dashed lines). In the right column of the figure, we show the metric quotients $\overline{MC_1(X)}, \overline{MC_2(X)}$ and $\overline{MC_3(X)}$ of the metric cylinder $MC(X)$. Though $MC(X)$ extends beyond $MC_3(X)$ with edges of zero weight they don't contribute to $\overline{MC(X)}$ and we have $\overline{MC(X)}=\overline{MC_3(X)}$. Note that $\overline{MC_3(X)}$ is a collection of trees (some of which are trivial).

\begin{figure}[H]
\centering
\scalebox{0.45}{
\begin{tabular}{@{}ccc@{}}
  \includegraphics[scale=1]{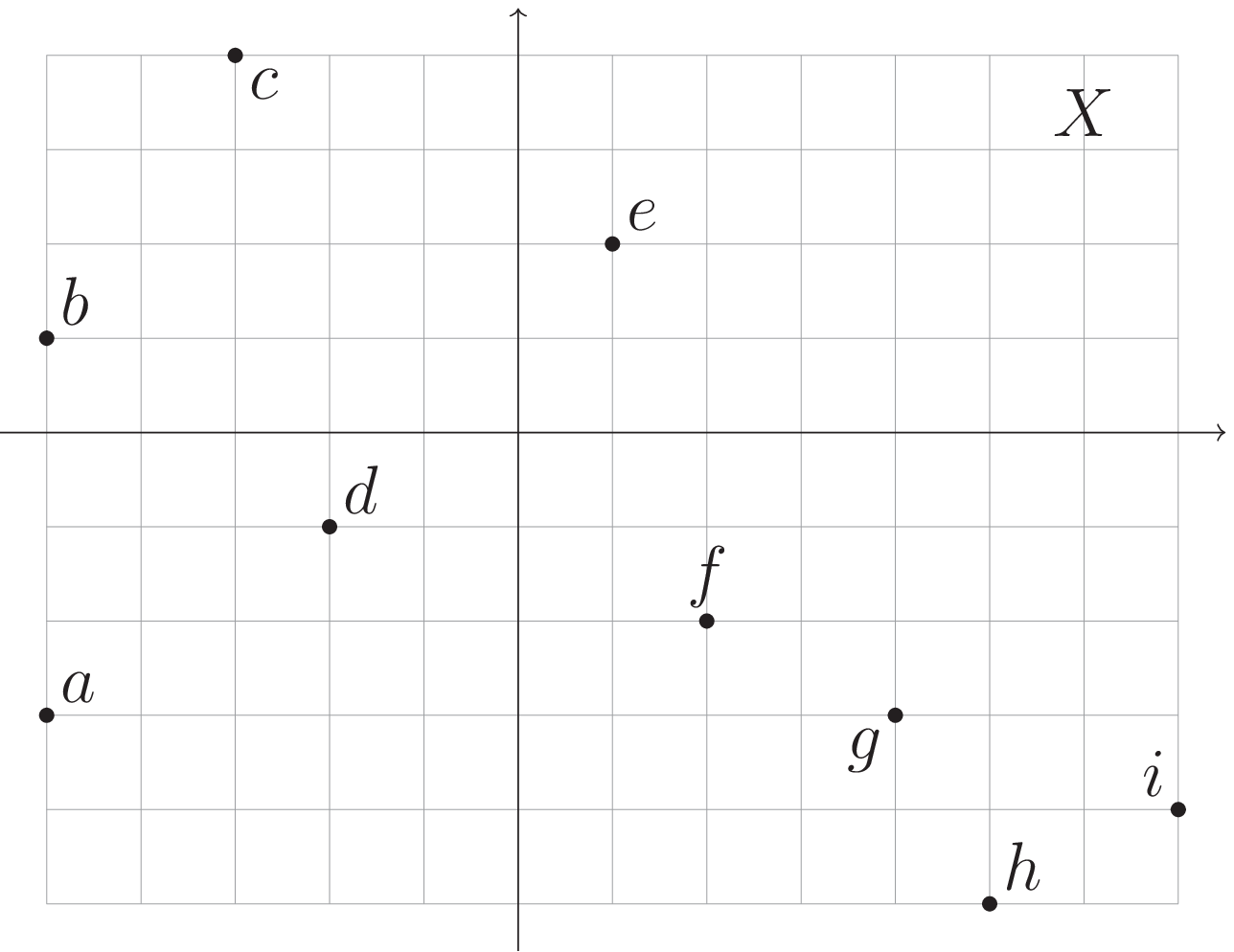} &\mbox{  }\hspace{1cm} &
  \includegraphics[scale=.1]{sekiz.eps} \\
      \includegraphics[scale=1]{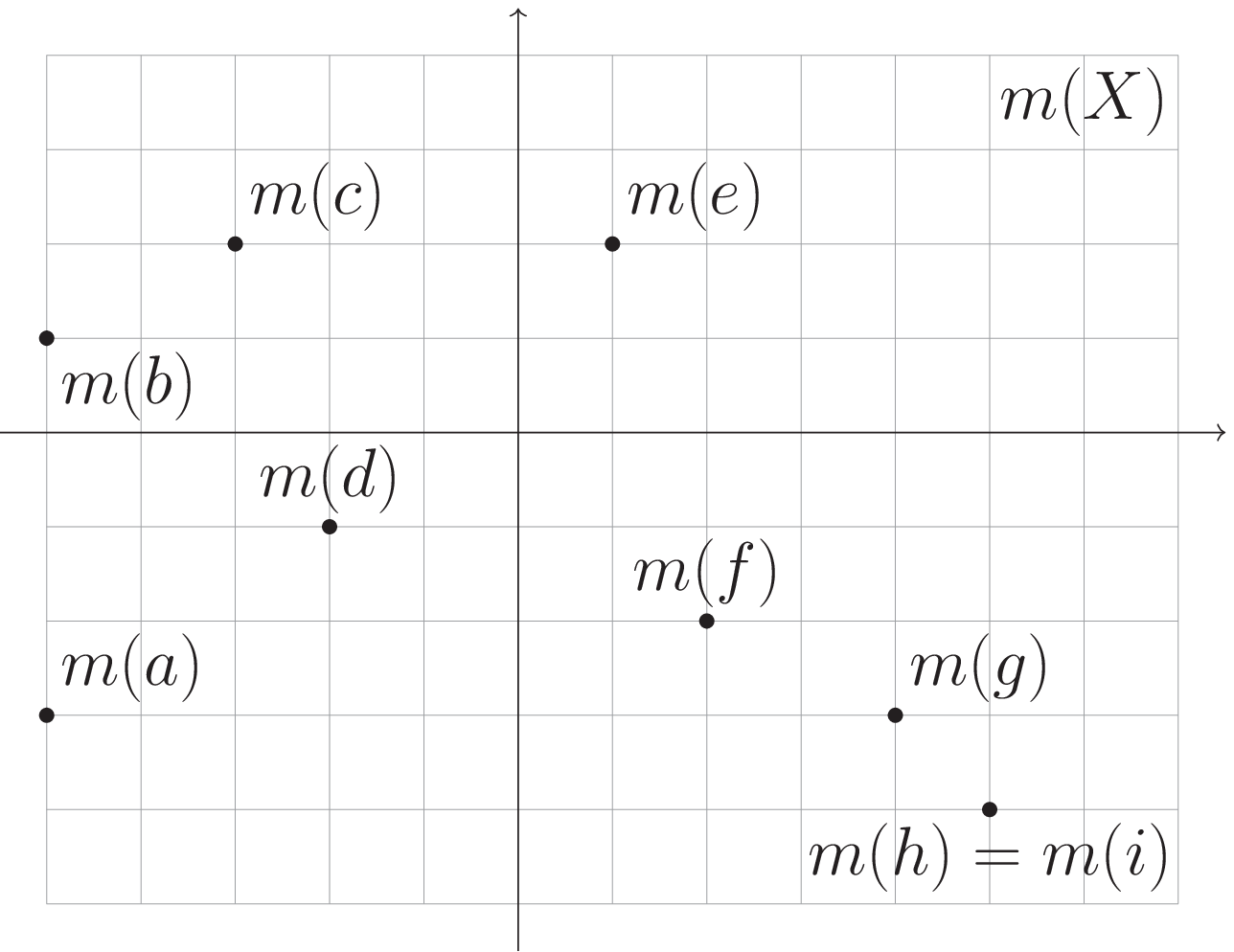} & \mbox{  }\hspace{1cm} &
   \includegraphics[scale=1]{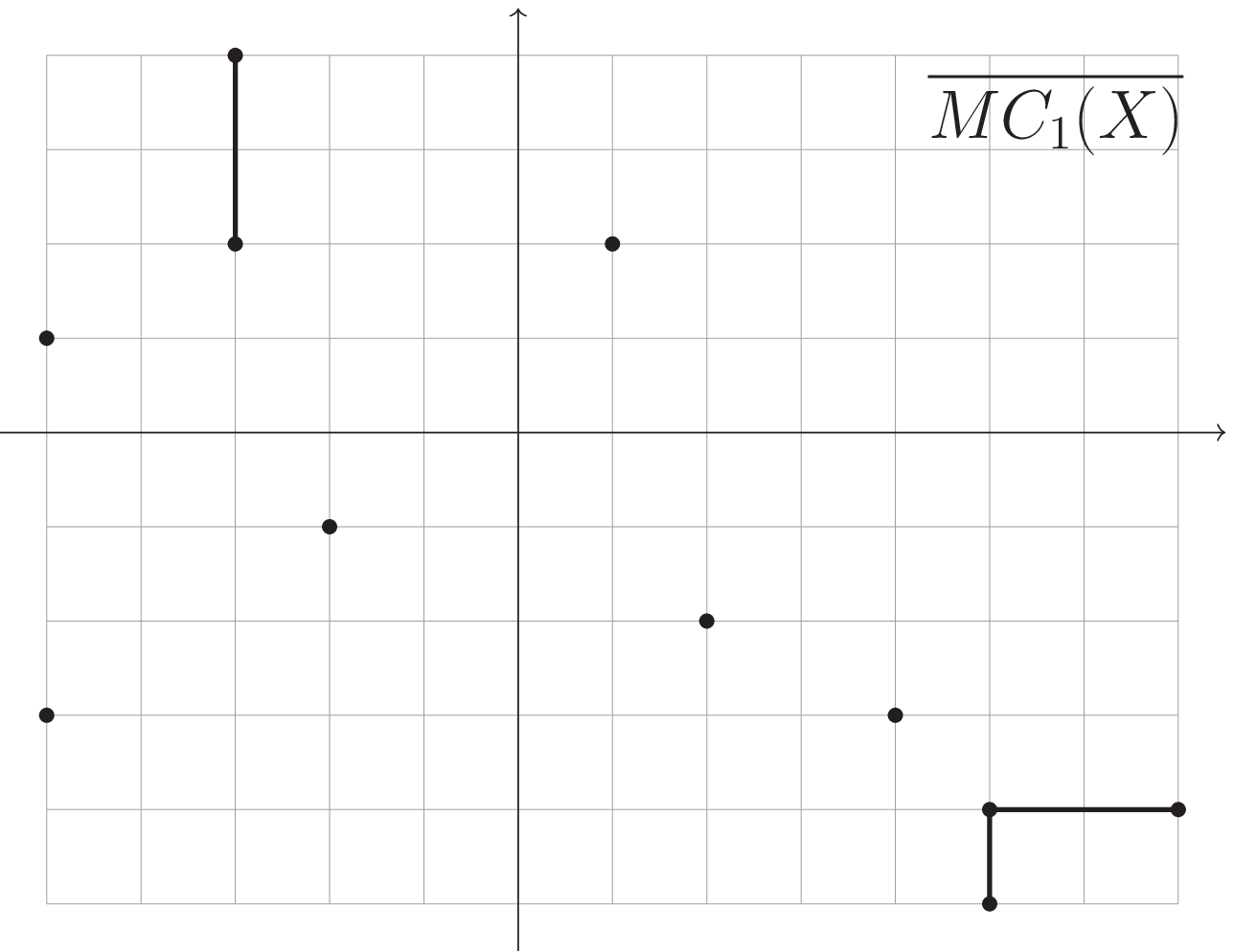} \\
       \includegraphics[scale=1]{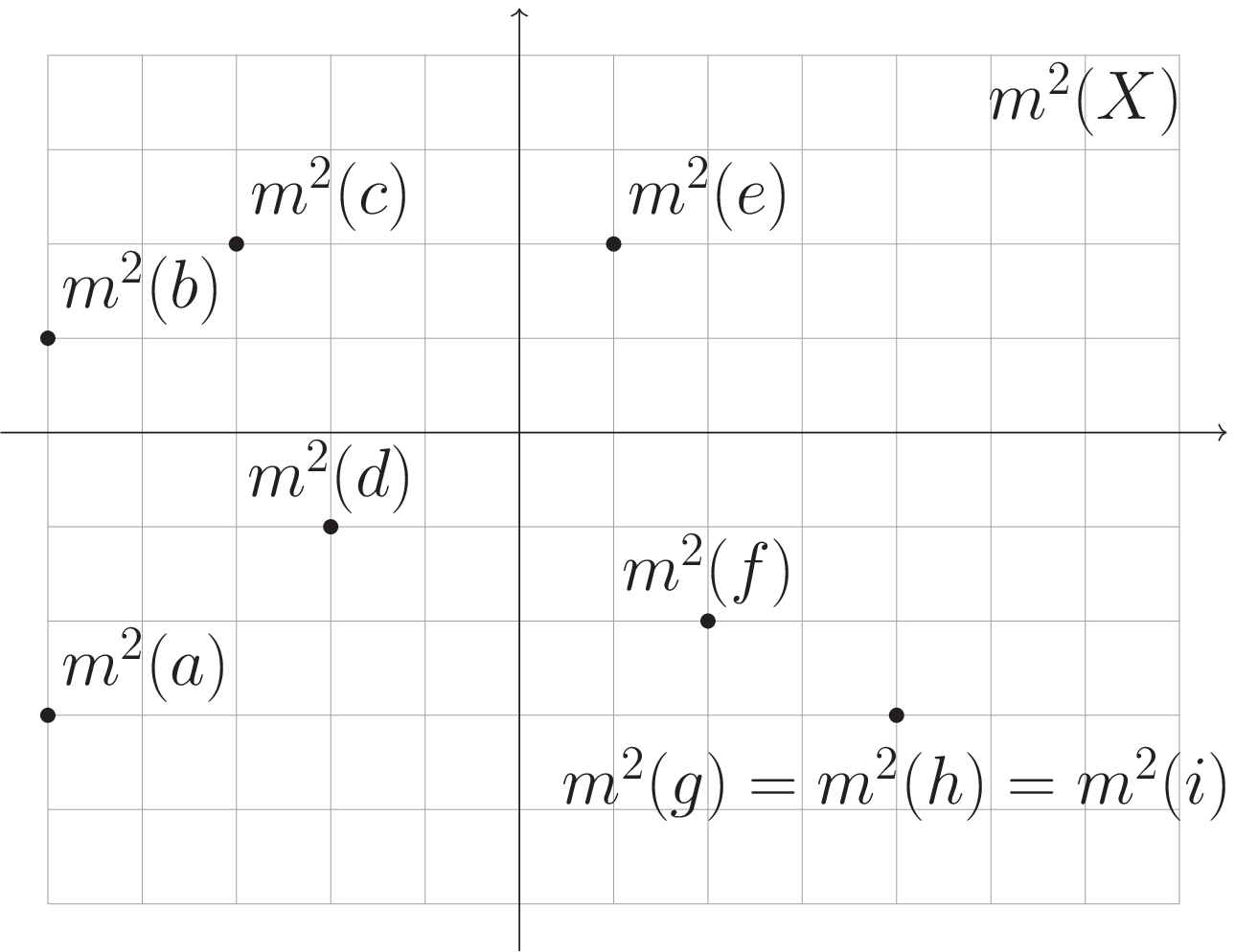} & \mbox{  } \hspace{1cm}&
           \includegraphics[scale=1]{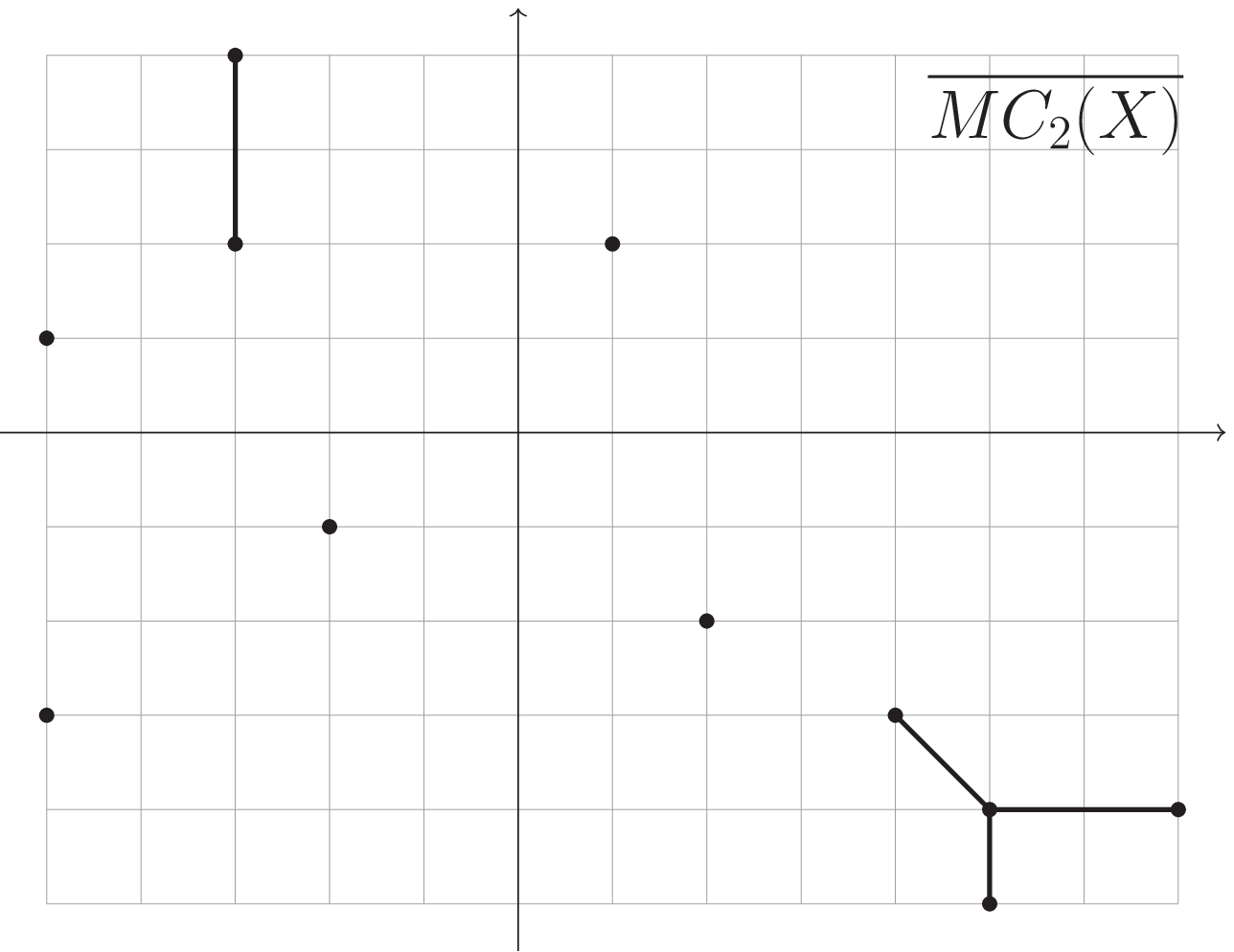} \\
         \includegraphics[scale=1]{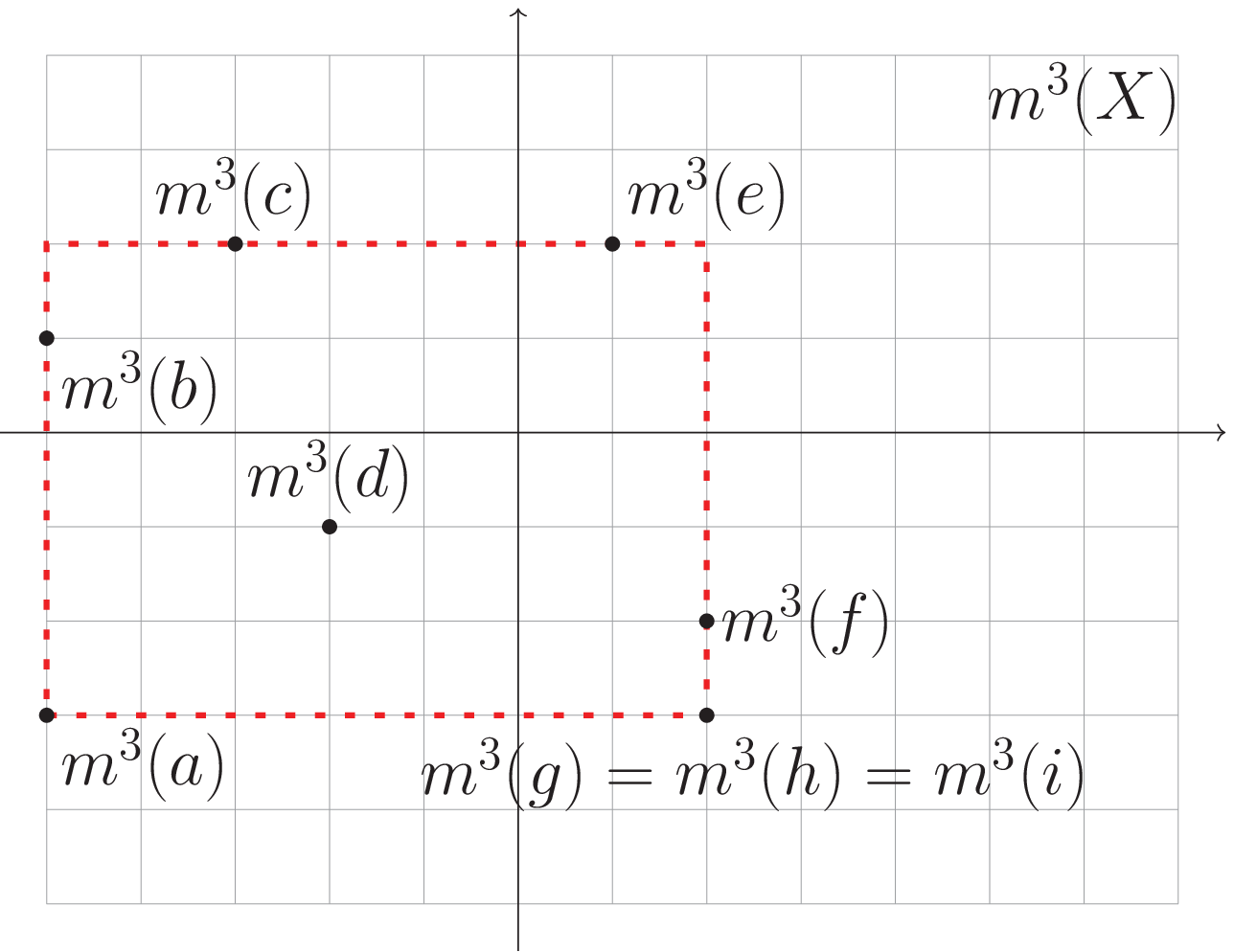} & \mbox{  } \hspace{1cm} &
           \includegraphics[scale=1]{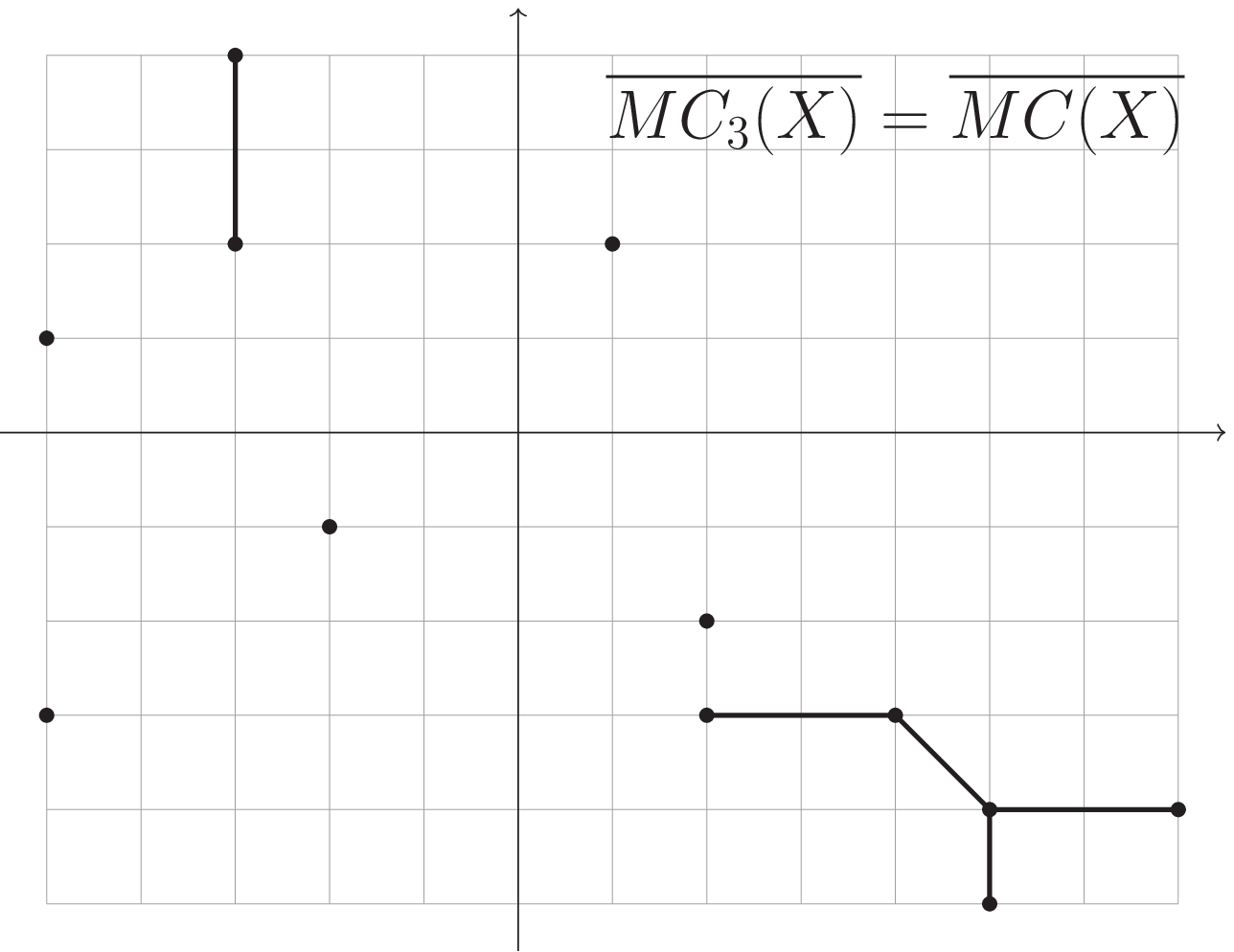} \\
  \end{tabular}}
  \caption{Successive metric centers and metric cylinders for the space $X$ of \\ Example~\ref{example1}.}\label{figex1}
\end{figure}
\end{example}

%%%%%
\section{The Tight Span of a Finite Subspace of the Manhattan Plane}
  In this section, we apply the main theorem of Turaev [\cite{Turaev2018},Theorem 9.1] to obtain the tight span of a finite subspace $X\subset \mathbb{R}^2_1$. Recall that, after successive trimming of $X$, a trim space is obtained (say after $N$ steps), so that the metric center sequence is stabilized at $m^N(X)\subset \mathbb{R}^2_1$. We then get by the main theorem of Turaev and our identification of successive trimmings by the successive metric center sets
  \begin{align*}
   T(X)&=\displaystyle \cap_{i\geq0}T(X_i) \cup \overline{C(X)}\\
   &=T(X_N) \cup \overline{C_N(X)}\\
   &=T(m^N(X)) \cup \overline{MC_N(X)}
  \end{align*}
  with $T(m^N(X)) \cap \overline{MC_N(X)}$ being the roots.

  Since we have described the $\overline{MC_N(X)}$ as an embedded graph in $\mathbb{R}^2_1$, which can be expressed as $\bigcup_{x\in X}\bigcup_{i=0}^{N-1}[m^i(x),m^{i+1}(x)]$, it remains to understand the term $T(m^N(X))$.

The stabilized (trim) $m^N(X)$ has either at least four points or a single point since a triple of points is not trim and likewise, a pair of points is not trim. If $m^N(X)$ is a singleton, then we get $T(X)=\overline{MC_N(X)}$, in which case $T(X)$ becomes a tree. Thus we can note the following property:
\begin{remark}
 The tight span of a finite subspace $X\subset\mathbb{R}^2_1$ is a tree if and only if after successive trimming of $X$ we get a singleton.

 (This property can easily be seen to be true for any finite metric space $X$.)
\end{remark}

In case the stabilized (trim) $m^N(X)$ has at least four points then, by Proposition \ref{theo0} the minimal rectangle $R_X$ has at least two points on every edge of it and we use the following theorem \cite{Kilic2016} to obtain the tight span of $m^N(X)$:
\begin{theorem}\label{theo6}
  Let $A\subseteq \mathbb{R}^2_1$ be a nonempty subspace. Let $B\subseteq \mathbb{R}^2_1$ be a closed, geodesically convex subspace containing $A$ and minimal with these properties. Then $B$ is isometric to the tight span $T(A)$ of $A$.
\end{theorem}

(This theorem is proven for the plane $\mathbb{R}^2_\infty$ with the maximum metric, but as $\mathbb{R}^2_\infty$ and $\mathbb{R}^2_1$ are isometric, it holds for $\mathbb{R}^2_1$ also.)

Now, we will develop an algorithm to create a subset containing $m^N(X)$ that is closed, geodesically convex and minimal with these properties to apply this theorem. Let $m^N(X)\subset \mathbb{R}^2_1$ be the initial set (with at least four points). Our algorithm to obtain $T(m^N(X))$ can be formulated as follows:
%%%%
(In the following, coordinates of a point $a\in \mathbb{R}^2_1$ is denoted by $a_x,a_y$, i.e. $a=(a_x,a_y)$.)
\begin{enumerate}
  \item [1.] Define $A:=m^N(X)$.
  \item [2.] Determine the minimal rectangle with edges parallel to the $x-$ and $y-$axis that contains all of the elements of $m^N(X)$ and call it $R$.
  \item [3.] Create a loop variable $P$ and define its initial value as $P:=R$.
  \item [4.] \textbf{Bottom Left Corner Procedure:} Choose the bottom left corner of $R$ and assign it to a variable $a=(a_x,a_y)$.
   \begin{enumerate}
     \item [(i)] If $a\in A$, then go to Step 5.
     \item [(ii)] If $a\notin A$, choose the closest point to $a$ with respect to the Manhattan metric from the set $A$ on the line $y=a_y$ and assign it to a variable $b=(b_x,b_y)$.
     \item [(iii)] Determine a point $d=(d_x,d_y)\in A$ (as a variable) such that $a_x\leq d_x\leq b_x$ and $a_y<d_y$ with smallest possible $d_y-a_y$.
     \item [(iv)] Remove $[a_x,b_x)\times [a_y,d_y)$ from $P$ and assign the remaining polygon to $P$. Add the point $(b_x,d_y)$ to the set $A$.
     \item [(v)] Choose $(a_x,d_y)$ as the new starting point and assign it to $a$. Return to Step 4(i).
   \end{enumerate}
  \item [5.] \textbf{Top Left Corner Procedure: }Choose the top left corner of $R$ and assign it to a variable $a=(a_x,a_y)$.
  \begin{enumerate}
     \item [(i)] If $a\in A$, then go to Step 6.
     \item [(ii)] If $a\notin A$, then choose the closest point to $a$ with respect to the Manhattan metric from the set $A$ on the line $y=a_y$ and assign it to a variable $b=(b_x,b_y)$.
     \item [(iii)] Determine a point $d=(d_x,d_y)\in A$ (as a variable) such that $a_x\leq d_x\leq b_x$ and $d_y<a_y$ with smallest possible $a_y-d_y$.
     \item [(iv)] Remove $[a_x,b_x)\times (d_y,a_y]$ from $P$ and assign the remaining polygon to $P$. Add the point $(b_x,d_y)$ to the set $A$.
     \item [(v)] Choose $(a_x,d_y)$ as the new starting point and assign it to $a$. Return to Step 5(i).
   \end{enumerate}
  \item [6.]  \textbf{Top Right Corner Procedure: }Choose the top right corner of $R$ and assign it to a variable $a=(a_x,a_y)$.
  \begin{enumerate}
    \item [(i)] If $a\in A$, then go to Step 7.
    \item [(ii)] If $a\notin A$, then choose the closest point to $a$ with respect to the Manhattan metric from the set $A$ on the line $y=a_y$ and assign it to a variable $b=(b_x,b_y)$.
    \item [(iii)] Determine a point $d=(d_x,d_y)\in A$ (as a variable) such that $b_x\leq d_x\leq a_x$ and $d_y<a_y$ with smallest possible $a_y-d_y$.
    \item [(iv)] Remove $(b_x,a_x]\times (d_y,a_y]$ from $P$ and assign the remaining polygon to $P$. Add the point $(b_x,d_y)$ to the set $A$.
    \item [(v)] Choose $(a_x,d_y)$ as the new starting point and assign it to $a$. Return to Step 6(i).
  \end{enumerate}
  \item [7.] \textbf{Bottom Right Corner Procedure:} Choose the bottom right corner of $R$ and assign it to a variable $a=(a_x,a_y)$.
   \begin{enumerate}
    \item [(i)] If $a\in A$, then the algorithm terminates.
    \item [(ii)] If $a\notin A$, then choose the closest point to $a$ with respect to the Manhattan metric from the set $A$ on the line $y=a_y$ and assign it to a variable $b=(b_x,b_y)$.
    \item [(iii)] Determine a point $d=(d_x,d_y)\in A$ (as a variable) such that $b_x\leq d_x\leq a_x$ and $a_y<d_y$ with smallest possible $d_y-a_y$.
    \item [(iv)] Remove $(b_x,a_x]\times [a_y,d_y)$ from $P$ and assign the remaining polygon to $P$. Add the point $(b_x,d_y)$ to the set $A$.
    \item [(v)] Choose $(a_x,d_y)$ as the new starting point and assign it to $a$. Return to Step 7(i).
  \end{enumerate}
\end{enumerate}

The source code of this algorithm along with the trimming procedure, which is written in the interactive geometry software Cinderella \cite{Cinderella1999} is given in the Appendix section.

%%%
\begin{example}\label{ex3}
Let $X$ be as in Figure \ref{fignew0}, i.e. $X=\{ (-5,-3),(-5,1),(-3,4),(-2,-1),$ \\$(1,2),(2,-2),(4,-3),(5,-5), (7,-4)\}$. We have obtained in Example \ref{example1} that $m^3(X)=\{(-5,-3),(-5,-1), (-3,2), (-2,-1), (1,2), (2,-3), (2,-2)\}$ is a trim space (see Figure \ref{figex2}(a)). Applying the above algorithm to this space, one will first obtain by the top left corner procedure Figure \ref{figex2}(b) and then by the top right corner procedure the tight span of $m^3(X)$ shown in Figure \ref{figex2}(c).

\begin{figure}[h]
\centering
\scalebox{.4}{
  \begin{tabular}{@{}ccc@{}}
  \includegraphics[scale=1]{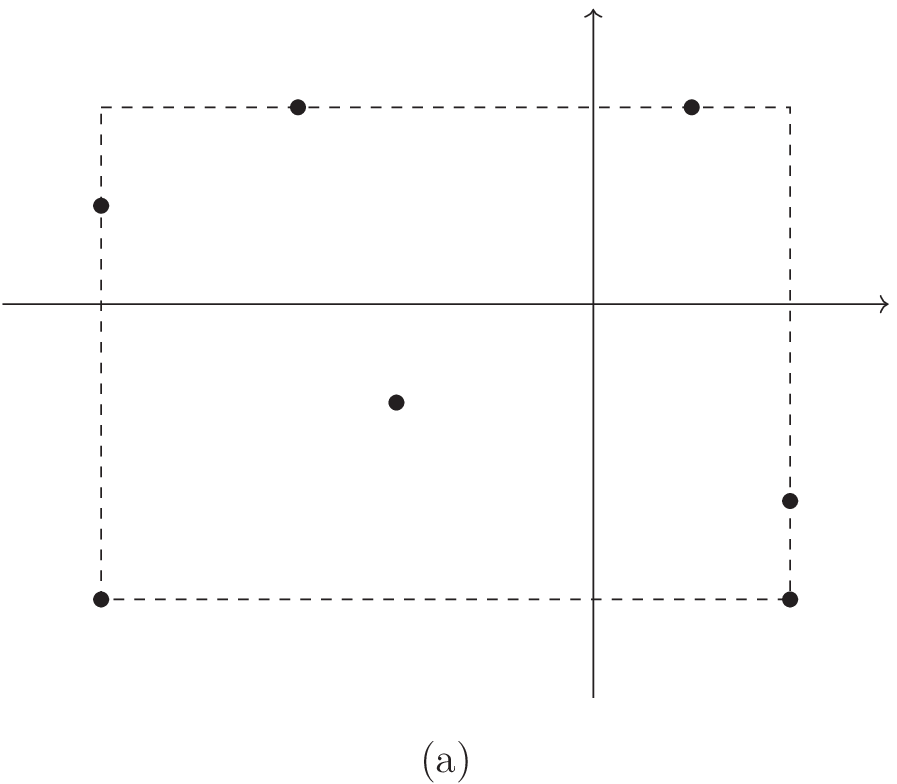} &
  \includegraphics[scale=1]{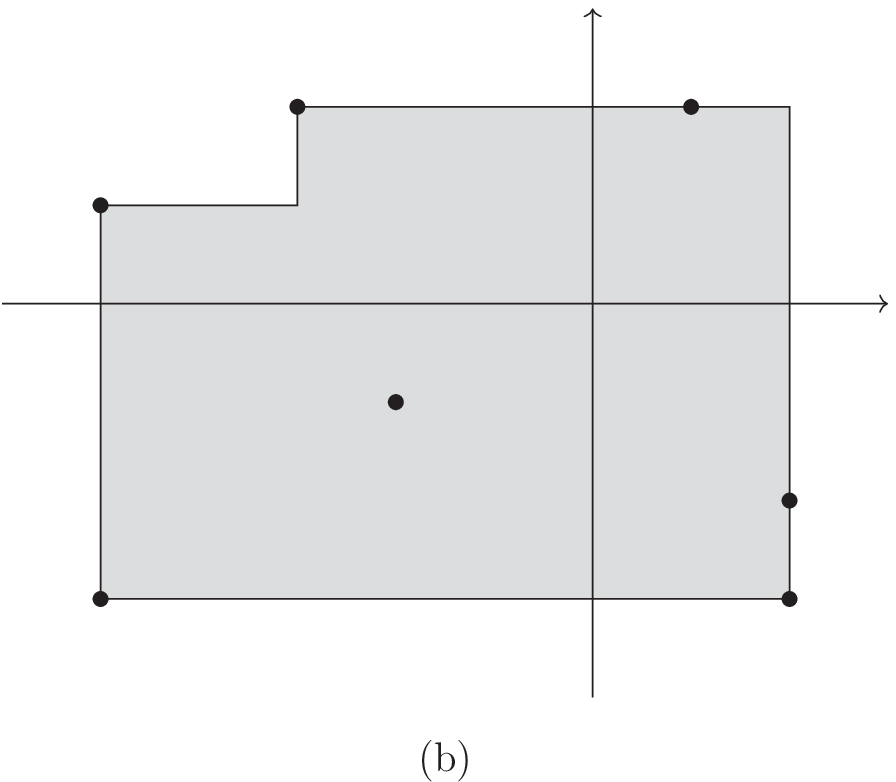}&
      \includegraphics[scale=1]{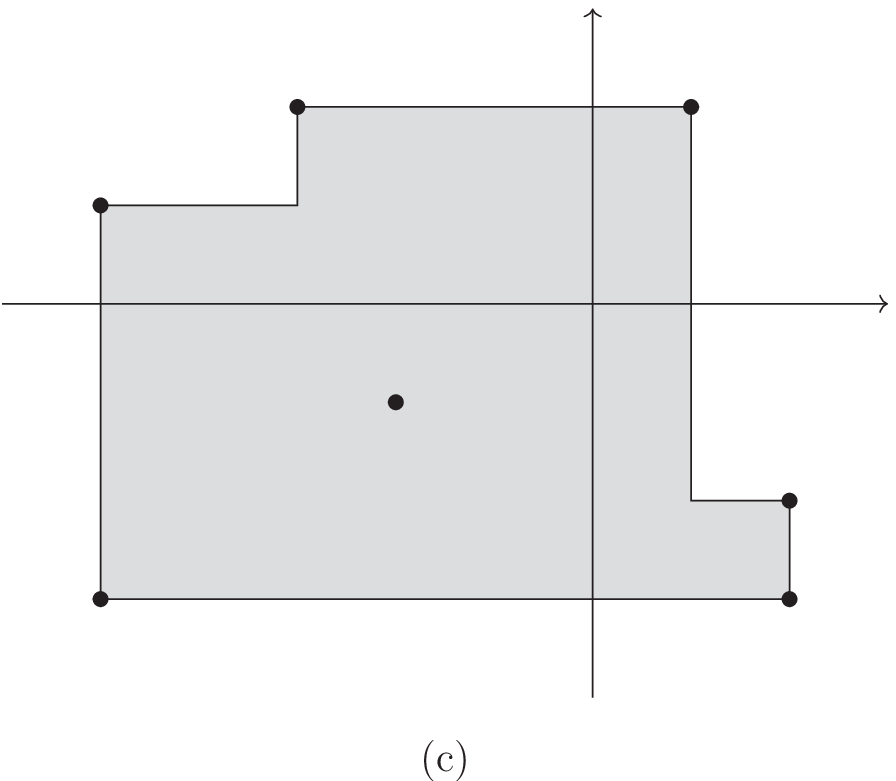}
  \end{tabular}}
  \caption{Application of the algorithm on Example \ref{example1}.}\label{figex2}
\end{figure}
\end{example}

After applying the algorithm to the metric space $m^N(X)$, denote the finally obtained set by $B$. The set $B$ contains $m^N(X)$ and it can be seen to be closed and geodesically convex. It is in fact also minimal since between two points with the same abcissa (or the same ordinate) there exists a unique geodesic (namely, the segment connecting them) and by this reason no point of $B$ can be discarded from the set. Having thus constructed the tight span of $m^N(X)$, we need only to take the sum of it with $\overline{MC_N(X)}$. Note that the intersection of these two sets consists of the points of $m^N(X)$ which are the ``roots" in the terminology of Turaev. The tips (leaves) of the attached trees belong to $X$ (some points of $X$ however might be intermediate points of the trees due to vanishing pendant lengths), and some trees might be constant (trivial) trees consisting of a single point, which is then a tip and a root simultaneously.

\begin{figure}[H]%m^2(X)
\centering
\scalebox{.57}{
\begin{tikzpicture} %[x={10.0pt},y={10.0pt}]
\foreach \i in {\xMin,...,\xMax} {
        \draw [very thin,gray!70] (\i,\yMin) -- (\i,\yMax);
    }
    \foreach \i in {\yMin,...,\yMax} {
        \draw [very thin,gray!70] (\xMin,\i) -- (\xMax,\i);
    }
  \draw [->, draw=gray!70] (-5.5,0) -- (8,0);  \draw [->,draw=gray!70] (0,-5.5) -- (0,4.5);
      \draw [fill=gray!30, draw=black, line width=0.6mm] (-5,1)--(-5,-3)--(2,-3)--(2,-2)--(1,-2)--(1,2)--(-3,2)--(-3,1)--cycle;
 \draw[black,fill=black](-5,1) circle (.4ex); \draw[black,fill=black](-3,4) circle (.4ex); \draw[black,fill=black](1,2) circle (.4ex);\draw[black,fill=black](-5,-3) circle (.4ex);
  \draw[black,fill=black](-2,-1) circle (.4ex); \draw[black,fill=black](2,-2) circle (.4ex); \draw[black,fill=black](4,-3) circle (.4ex);
 \draw[black,fill=black](5,-5) circle (.4ex);\draw[black,fill=black](7,-4) circle (.4ex);

 \draw[line width=0.7mm](-3,4) -- (-3,2);  \draw[line width=0.7mm](5,-4) -- (7,-4);  \draw[line width=0.7mm](5,-5) -- (5,-4); \draw[line width=0.7mm](4,-3) -- (5,-4); \draw[line width=0.7mm](4,-3) -- (2,-3);
\end{tikzpicture}}
\caption{The tight span of $X=\{ (-5,-3),(-5,1),(-3,4),(-2,-1), (1,2),$ $(2,-2),(4,-3),(5,-5), (7,-4)\}$.}\label{figex3}
\end{figure}
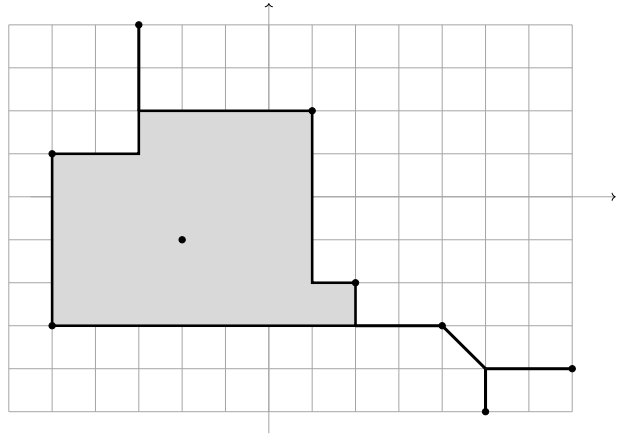

We exemplify this last step of constructing the tight span of a finite subset of $\mathbb{R}^2_1$ on our ongoing example subspace $X$ (of Example \ref{example1}). We had constructed the $\overline{MC(X)}$ in Example \ref{example1} (see Figure \ref{figex1}) and we have above constructed the tight span of $m^3(X)$ (see Figure \ref{figex2}(c)) so that we can now take their union to obtain the tight span of $X$ (see Figure \ref{figex3}).

Finally, we want to conclude with two more examples:
\begin{example}\label{ex4}
 Let $X=\{(-5,2),(-1,0),(-1,3),(3,2),(4,-3),(5,-1),(6,-4),$ \\ $(8,-5)\}$. The tight span of $X$ is shown in Figure \ref{fig30}. The metric center sequence stabilizes at $m^2(X)=\{(-1,0),(-1,2),(3,2),(4,-3),(5,-3),(5,-1)\}$.
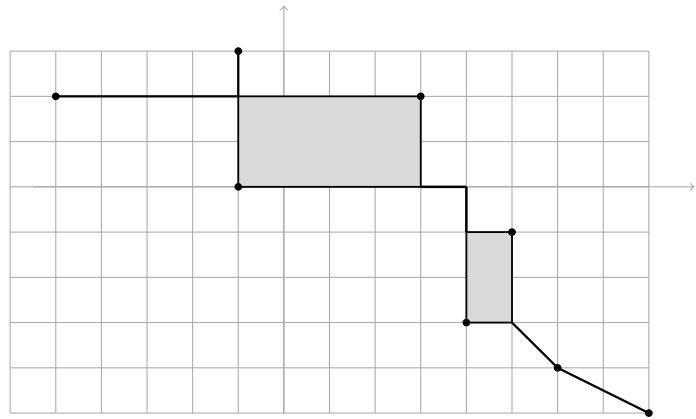
\begin{figure}[H]
\centering
\scalebox{.6}{
\begin{tikzpicture} %[x={10.0pt},y={10.0pt}]
\foreach \i in {-6,...,8} {
        \draw [very thin,gray!70] (\i,-5) -- (\i,3);
    }
    \foreach \i in {-5,...,3} {
        \draw [very thin,gray!70] (-6,\i) -- (8,\i);
    }
 \draw [gray!70,->] (-5.5,0) -- (9,0);  \draw [gray!70,->] (0,-5) -- (0,4);
    \draw [fill=gray!30, draw=black,line width=0.4mm] (-1,0)--(3,0)--(3,2)--(-1,2)--cycle;
    \draw [fill=gray!30, draw=black,line width=0.4mm] (4,-3)--(5,-3)--(5,-1)--(4,-1)--cycle;

 \draw[black,fill=black](-5,2) circle (.4ex);  \draw[black,fill=black](-1,0) circle (.4ex);  \draw[black,fill=black](-1,3) circle (.4ex);  \draw[black,fill=black](3,2) circle (.4ex);  \draw[black,fill=black](4,-3) circle (.4ex); \draw[black,fill=black](5,-1) circle (.4ex);  \draw[black,fill=black](6,-4) circle (.4ex);  \draw[black,fill=black](8,-5) circle (.4ex);

  \draw[line width=0.5mm](-5,2)--(-1,2);  \draw[line width=0.5mm](-1,2)--(-1,3);  \draw[line width=0.5mm](5,-3)--(6,-4);  \draw[line width=0.5mm](6,-4)--(8,-5);

 \draw[line width=0.55mm](3,0)--(4,0);  \draw[line width=0.55mm](4,-1)--(4,0);

\end{tikzpicture}}
\caption{Tight span of $X$ of Example \ref{ex4}.}\label{fig30}
\end{figure}
\end{example}

\begin{example}\label{ex5}
 Let $X=\{(-15,7),(-13,10),(-12,5),(-9,-3),(-9,0),(-7,4),$ \\ $(-5,-5),(-5,4),(-3,-7),(-2,-2),(-1,-9),(-1,2),(1,4),(2,-5),(4,-7),$ \\ $(4,-4),(7,-6)\}.$ The tight span of $X$ is shown in Figure \ref{fig27}. The metric center sequence stabilizes at $m^3(X)=\{(-9,-3),(-9,0),(-9,4),(-7,4),(-5,-5),$ \\ $(-5,4),(-3,-7),(-2,-2),(-1,-7),(-1,2),(1,4),(2,-5),(4,-7),(4,-6),$ \\ $(4,-4)\}$.
\begin{figure}[H]
\centering
\scalebox{.35}{
\begin{tikzpicture} %[x={10.0pt},y={10.0pt}]
\foreach \i in {-16,...,8} {
        \draw [very thin,gray!70] (\i,-10) -- (\i,10);
    }
    \foreach \i in {-10,...,10} {
        \draw [very thin,gray!70] (-16,\i) -- (8,\i);
    }

  \draw [gray,->] (-15.5,0) -- (9.5,0);  \draw [gray,->] (0,-10.5) -- (0,11.5);
    \draw [fill=gray!30, draw=black,line width=0.7mm] (-9,4)--(-9,-3)--(-5,-3)--(-5,-5)--(-3,-5)--(-3,-7)--(4,-7)--(4,-4)--(1,-4)--(1,4)--cycle;

     \draw[black,fill=black](-1,2) circle (.5ex); \draw[black,fill=black](1,4) circle (.5ex); \draw[black,fill=black](-5,4) circle (.5ex);\draw[black,fill=black](-7,4) circle (.5ex);
  \draw[black,fill=black](-9,0) circle (.5ex); \draw[black,fill=black](-9,-3) circle (.5ex); \draw[black,fill=black](-2,-2) circle (.5ex);\draw[black,fill=black](4,-4) circle (.5ex);
  \draw[black,fill=black](2,-5) circle (.5ex); \draw[black,fill=black](-5,-5) circle (.5ex); \draw[black,fill=black](-3,-7) circle (.5ex);\draw[black,fill=black](4,-7) circle (.5ex);
   \draw[black,fill=black](-1,-9) circle (.5ex); \draw[black,fill=black](7,-6) circle (.5ex); \draw[black,fill=black](-12,5) circle (.5ex);\draw[black,fill=black](-13,10) circle (.5ex);
 \draw[black,fill=black](-15,7) circle (.5ex);

  \draw[line width=0.8mm](-13,7)--(-15,7);  \draw[line width=0.8mm](-13,7)--(-13,10); \draw[line width=0.8mm](-13,7)--(-12,5); \draw[line width=0.8mm](-12,5)--(-9,4);  \draw[line width=0.8mm](-1,-9)--(-1,-7); \draw[line width=0.8mm](7,-6)--(4,-6);
\end{tikzpicture}}
\caption{Tight span of $X$ of Example \ref{ex5}.}\label{fig27}
\end{figure}
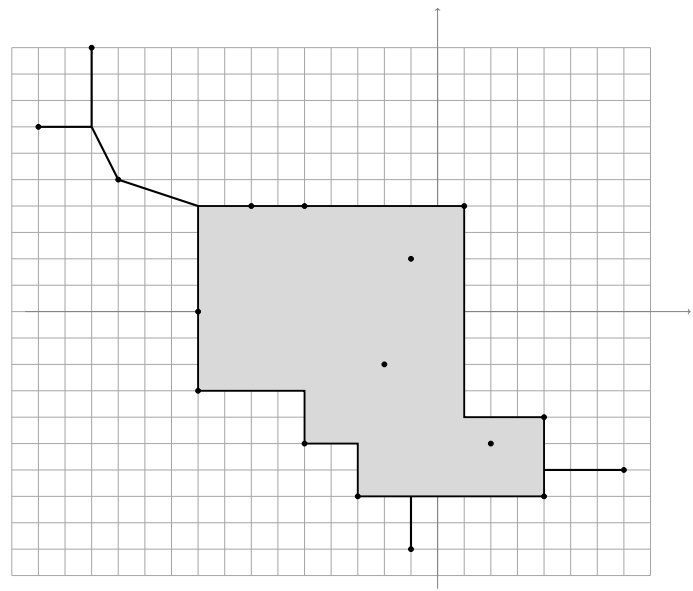
\end{example}
%%%%%%

\section{Appendix}
To use this code, copy and paste it page by page into ``Draw" under Cinderella's Script Editor.

\begin{verbatim}
ptsall=allpoints();
trimseq=[ptsall];
//***************************************************
k:=1;
pts=trimseq_k;
repeat(length(pts),i,pts_i=gauss(complex(pts_i)));
infty=1000000;
xmax1=-infty; xmax2=-infty-1; ymax1=-infty; ymax2=-infty-1;
xmin1=infty;  xmin2=infty-1;  ymin1=infty;  ymin2=infty-1;
while((length(pts)>2 & ((xmax1 ~!= xmax2) % (xmin1 ~!= xmin2)
% (ymax1 ~!= ymax2) %(ymin1 ~!= ymin2))),
    xmax1=-infty; xmax2=-infty-1; ymax1=-infty; ymax2=-infty-1;
    xmin1=infty;  xmin2=infty-1;  ymin1=infty;  ymin2=infty-1;
repeat (length(pts),a,
if ((pts_a)_1 >= xmax1,
   xmax2=xmax1;
   xmax1=(pts_a)_1
   ,if ((pts_a)_1 > xmax2, xmax2=(pts_a)_1)
 );
if ((pts_a)_1 <= xmin1,
   xmin2=xmin1;
   xmin1=(pts_a)_1
 ,if ((pts_a)_1 < xmin2, xmin2=(pts_a)_1)
 );
if ((pts_a)_2 >= ymax1,
   ymax2=ymax1;
   ymax1=(pts_a)_2
 ,if ((pts_a)_2 > ymax2, ymax2=(pts_a)_2)
 );
if ((pts_a)_2 <= ymin1,
   ymin2=ymin1;
   ymin1=(pts_a)_2
 ,if ((pts_a)_2 < ymin2, ymin2=(pts_a)_2)
 );
);
trimpts=pts;
repeat(length(pts),i,
if ((trimpts_i)_1 <= xmin1, (trimpts_i)_1=xmin2;);
if ((trimpts_i)_1 >= xmax1, (trimpts_i)_1=xmax2;);
if ((trimpts_i)_2 <= ymin1, (trimpts_i)_2=ymin2;);
if ((trimpts_i)_2 >= ymax1, (trimpts_i)_2=ymax2;);
draw((pts)_i,(trimpts)_i);
);
trimseq=append(trimseq,trimpts);
k=k+1;
pts=set(trimseq_k);
repeat(length(pts),i,draw(pts_i));
);
if (length(pts)==2,
trimpts=[[((pts_1)_1+(pts_2)_1)/2,(trimpts_1)_2=((pts_1)_2+(pts_2)_2)/2]];
draw(pts_1,trimpts_1); draw(pts_2,trimpts_1);
draw(trimpts_1);
trimseq=append(trimseq,trimpts);
);
trimstep=length(trimseq);
if (length(trimseq_trimstep)>=4,
tightspan=[[xmin1,ymin1],[xmin1,ymax1],[xmax1,ymax1],[xmax1,ymin1]];
pts=trimseq_trimstep;

// BOTTOM-LEFT CORNER PROCEDURE
refpoint=[xmin1,ymin1];
while(!contains(pts,refpoint),
hat=select(pts,#.y==refpoint.y);
xfirst=[infty,refpoint.y];
repeat(length(hat),i,if ((xfirst.x > (hat_i).x),  xfirst=hat_i); );
strip=select(pts,((#.y > xfirst.y) & (#.x < xfirst.x)));
yfirst=[xmin1,ymax1];
repeat(length(strip),i, if ((yfirst.y > (strip_i).y), yfirst=strip_i); );
init=[]; final=[]; mid=[xfirst,[xfirst.x,yfirst.y],[xmin1.x,yfirst.y]];
draw(mid_1,mid_2);
draw(mid_2,[yfirst.x,yfirst.y]);
control=0;
repeat (length(tightspan),k,
if (tightspan_k==refpoint,control=1; );
if (control==0, init=init ++[tightspan_k]);
if ((control==1)&(tightspan_k!=refpoint),final=final++[tightspan_k]);
);
tightspan=concat(concat(init,mid),final);
pts=append(pts,[xfirst.x,yfirst.y]);
refpoint=mid_3;
);
// TOP-LEFT CORNER PROCEDURE
refpoint=[xmin1,ymax1];
while(!contains(pts,refpoint),
hat=select(pts,#.y==refpoint.y);
xfirst=[infty,refpoint.y];
repeat(length(hat),i,if ((xfirst.x > (hat_i).x),  xfirst=hat_i); );
strip=select(pts,((#.y < xfirst.y)   & (#.x<=xfirst.x) ));
yfirst=[xmin1,ymin1];
repeat(length(strip),i, if ((yfirst.y < (strip_i).y), yfirst=strip_i); );
init=[]; final=[];  mid=[[xmin1.x,yfirst.y],[xfirst.x,yfirst.y],xfirst];
draw(mid_3,mid_2);
draw(mid_2,[yfirst.x,yfirst.y]);
control=0;
repeat (length(tightspan),k,
if (tightspan_k==refpoint,control=1; );
if (control==0, init=init++[tightspan_k]);
if ((control==1)&(tightspan_k!=refpoint),final=final++[tightspan_k]);
);

tightspan=concat(concat(init,mid),final);
pts=append(pts,[xfirst.x,yfirst.y]);
refpoint=mid_1;
);
// TOP-RIGHT CORNER PROCEDURE
refpoint=[xmax1,ymax1];
while(!contains(pts,refpoint),
hat=select(pts,#.y==refpoint.y);
xfirst=[-infty,refpoint.y];
repeat(length(hat),i,if ((xfirst.x < (hat_i).x),  xfirst=hat_i); );
strip=select(pts,((#.y < xfirst.y)   & (#.x >= xfirst.x) ));
yfirst=[xmax1,ymin1];
repeat(length(strip),i, if ((yfirst.y < (strip_i).y), yfirst=strip_i); );
init=[];
final=[];
mid=[xfirst,[xfirst.x,yfirst.y],[xmax1.x,yfirst.y]];
draw(mid_1,mid_2);
draw(mid_2,[yfirst.x,yfirst.y]);
control=0;
repeat (length(tightspan),k,
if (tightspan_k==refpoint,control=1; );
if (control==0, init=init++[tightspan_k]);
if ((control==1)&(tightspan_k!=refpoint),final=final++[tightspan_k]);
);

tightspan=concat(concat(init,mid),final);
pts=append(pts,[xfirst.x,yfirst.y]);
refpoint=mid_3;
);
// BOTTOM-RIGHT CORNER PROCEDURE
refpoint=[xmax1,ymin1];
while(!contains(pts,refpoint),
hat=select(pts,#.y==refpoint.y);
xfirst=[-infty,refpoint.y];
repeat(length(hat),i,if ((xfirst.x < (hat_i).x),  xfirst=hat_i); );
strip=select(pts,((#.y > xfirst.y)   & (#.x >= xfirst.x) ));
yfirst=[xmax1,ymax1];
repeat(length(strip),i, if ((yfirst.y > (strip_i).y), yfirst=strip_i); );
init=[];
final=[];
mid=[[xmax1.x,yfirst.y],[xfirst.x,yfirst.y],xfirst];
draw(mid_3,mid_2);
draw(mid_2,[yfirst.x,yfirst.y]);
control=0;
repeat (length(tightspan),k,
if (tightspan_k==refpoint,control=1; );
if (control==0, init=init++[tightspan_k]);
if ((control==1)&(tightspan_k!=refpoint),final=final++[tightspan_k]);
);
tightspan=concat(concat(init,mid),final);
pts=append(pts,[xfirst.x,yfirst.y]);
refpoint=mid_1;
);

tightspan=append(tightspan,tightspan_1);
repeat(length(tightspan)-1,i,
draw(tightspan_i,tightspan_(i+1));
draw(tightspan_i);
);
alpha(0.3);
fillpoly(tightspan,color->hue(0.52));
);


\end{verbatim}

%%%%%%%%%%%%%%%%%%%%%%%%%%


\begin{thebibliography}{9}

\bibitem{Bilge2017} Bilge, A. H., Çelik, D., Koçak, Ş., An equivalence class decomposition of finite metric spaces via Gromov products, \textit{Discrete Mathematics}, 340(8), 1928–1932, (2017).

\bibitem{Kilic2016} Kılıç, M. and Koçak, Ş., Tight span of subsets of the plane with the maximum metric. \textit{Advances in Mathematics}, 301, 693-710, (2016).

\bibitem{Kilic2021} Kılıç, M., Koçak, Ş., Özdemir, Y., An algorithm for the construction of the tight span of finite subsets of the Manhattan plane. \textit{Computational Geometry}, Vol. 95, 101741, (2021).

\bibitem{Cinderella1999} Richter-Gebert J., Kortenkamp U., The interactive geometry software {C}inderella, Springer-Verlag, Berlin, (1999).

\bibitem{Turaev2016} Turaev, V., Trimming of finite metric spaces. \textit{arXiv:1612.06760v1}, (2016).

\bibitem{Turaev2018} Turaev, V., Trimming of metric spaces and the tight span. \textit{Discrete Mathematics}, 341(10), 2925-2937, (2018).

\end{thebibliography}
\end{document}